\documentclass{amsart}

\usepackage{graphics}
\usepackage{amsmath,amsfonts}
\usepackage{amstext} 
\usepackage{amsbsy,amsxtra}
\usepackage{enumerate}
\usepackage{amssymb}
\usepackage{latexsym}
\usepackage{graphicx}

 \usepackage[dvips]{color}

 \newtheorem{theorem}{Theorem}[section]
 \newtheorem{proposition}[theorem]{Proposition}
 \newtheorem{lemma}[theorem]{Lemma}
 \newtheorem{corollary}[theorem]{Corollary}

 \theoremstyle{definition}
 \newtheorem{definition}[theorem]{Definition}
 \newtheorem{example}[theorem]{Example}
 \newtheorem{remark}[theorem]{Remark}

\numberwithin{theorem}{section} \numberwithin{equation}{section}

\newcommand{\nR}{{\mathbb R}}

\newcommand{\cJ}{{\mathcal J}}
\newcommand{\cF}{{\mathcal F}}

\newcommand{\wh}[1]{\widehat{#1}}

\newcommand{\bL}{{\mathbf L}}
\newcommand{\bK}{{\mathbf K}}

 \DeclareMathOperator{\co}{co}
 \DeclareMathOperator{\conv}{conv}
 \DeclareMathOperator{\aff}{aff}

\begingroup
\def\2#1{\ifnum#1<10 0\fi\the#1}
\xdef\isodayandtime%
{\the\year-\2\month-\2\day\space\2{\count0}:\2{\count2}}

\endgroup
\newcommand{\rr}{\mathbb{R}}

\newcommand{\su}{\sum_{j=1}^{m}}
\newcommand{\jom}{j=1,\ldots,m}

\begin{document}

\title[Convex Operator]%
{On a Convex Operator for Finite Sets}

\author{Branko \'{C}urgus}
\address{Department of Mathematics, Western Washington University,
Bellingham, Washington 98226, USA}
 \email{curgus@cc.wwu.edu}

\author{Krzysztof Ko{\l}odziejczyk}
 \address{Institute of Mathematics and Computer Science,
 Wroc{\l}aw University of Technology, 50-370 Wroc{\l}aw, Poland}
 \email{Krzysztof.Kolodziejczyk@pwr.wroc.pl}



\begin{abstract}
Let $S$ be a finite set with $n$ elements in a real linear space.
Let $\cJ_S$ be a set of $n$ intervals in $\nR$.  We introduce a
convex operator $\co(S,\cJ_S)$ which generalizes the familiar
concepts of the convex hull $\conv S$ and the affine hull $\aff S$
of $S$.  We establish basic properties of this operator.  It is
proved that each homothet of $\conv S$ that is contained in $\aff S$
can be obtained using this operator.  A variety of convex subsets of
$\aff S$ can also be obtained. For example, this operator assigns a
regular dodecagon to the $4$-element set consisting of the vertices
and the orthocenter of an equilateral triangle.  For $\cJ_S$ which
consists of bounded intervals, we give the upper bound for the
number of vertices of the polytope $\co(S,\cJ_S)$.
\end{abstract}

\maketitle


\section{Introduction}
Let $d$ be a positive integer and let $\bL$ be a real linear space
of dimension $d$.  Let $S=\{x_1,\dots,x_m\}$ be a finite set of
distinct points in $\bL$.  The familiar convex sets associated with
$S$ are the convex hull of $S$
\begin{align} \label{edconv}
\conv S & = \Biggl\{ \sum\limits_{j=1}^{m} \xi_j\, x_j \, : \, x_j
\in S, \  \xi_j \ge 0, \ \sum\limits_{j=1}^m \xi_j = 1
\Biggr\} \\
 \intertext{and the affine hull of $S$}
 \nonumber 
\aff S & = \Biggl\{\sum\limits_{j=1}^{m} \xi_j\, x_j \, : \,
 x_j \in S, \ \xi_j \in \nR, \
 \sum\limits_{j=1}^m \xi_j = 1\Biggr\}.
\end{align}

The set in \eqref{edconv} will not change if the conditions $\xi_j
\geq 0$ are replaced by the conditions $\xi_j \in [0,1], \,
j=1,\ldots,m$. This leads us to ask the following natural question:
How will the set on the right-hand side of \eqref{edconv} change
when the conditions $\xi_j \ge 0$ are replaced by the conditions
$\xi_j \in I_j, \, j = 1, \ldots, m$, where $I_j$ are arbitrary
nonempty intervals in $\nR$?  An immediate and obvious answer is
that the resulting set will always be a subset of $\aff S$. In this
article we explore this question further. That is, we study the
subsets of $\bL$ introduced by the following definition.

\begin{definition} \label{dco}
Let $S = \{x_1,\ldots,x_m\}$ be a finite set of distinct points in a
linear space $\bL$. Let $\cJ_S = \{ I_1, \ldots, I_m\}$ be a family
of nonempty intervals in $\nR$ (some of which can be degenerated to
a singleton) such that the interval $I_j$ is associated with the
point $x_j$ for each $j\in \{1,\ldots,m\}$. Set
\begin{equation*}
\co(S,\cJ_S) := \Biggl\{ \sum\limits_{j=1}^{m} \xi_j\, x_j \, : \,
x_j \in S, \  \xi_j \in I_j, \ \sum \limits_{j=1}^m \xi_j = 1
\Biggr\}.
\end{equation*}
This set we call a {\em convex interval hull} of $S$.
\end{definition}

It is clear that the convex interval hull $\co(S,\cJ_S)$ coincides
with $\conv S$ when all the intervals in $\cJ_S$ are equal to
$[0,1]$ and $\co(S,\cJ_S)$ coincides with $\aff S$ when all the
intervals in $\cJ_S$ are equal to $\nR$.  In this sense
$\co(S,\cJ_S)$ generalizes these two well-known concepts.

Our primary interest in this article is to explore the family of all
convex interval hulls $\co(S,\cJ_S)$ which are bounded.  Examples in
Section~\ref{sdbp} show that a variety of convex sets appear in such
families even if $S$ is fixed.  It is quite striking that when $S$
is the set of only $4$ points: the vertices and the orthocenter of
an equilateral triangle, for example
\begin{align*}
S & = \bigl\{\bigl(-1,0\bigr), \bigl(0,1\bigr),
\bigl(0,\sqrt{3}\,\bigr), \bigl( 0,1/\sqrt{3}\, \bigr) \bigr\},  \\
\intertext{and when} \cJ_S & = \bigl\{ [0,1], [0,1],
[0,1],\bigl[1-\sqrt{3}, -2 + \sqrt{3}\bigr] \bigr\},
\end{align*}
then $\co(S,\cJ_S)$ is a regular dodecagon; see Figure~\ref{f12}.

An inverse problem in this setting is as follows: For a given convex
set $K$ find a finite set $S$ with minimal cardinality and a family
$\cJ_S$ of intervals such that $K = \co(S,\cJ_S)$.
Examples~\ref{ex3} and \ref{exrhdo} suggest solutions of the inverse
problem for $K$ equal to a regular dodecagon and for $K$ equal to a
rhombic dodecahedron. This inverse problem and unbounded convex
interval hulls will be considered elsewhere.

Let $m$ be a positive integer.  For a fixed family $\cJ$ of $m$
nonempty intervals in $\nR$ our operator $S \mapsto \co(S,\cJ)$ is a
set-valued function defined on finite subsets of $\bL$ with $m$
elements. Recall that many set-valued functions $f$ considered in
convexity theory are described in the following way:
\begin{equation} \label{eqff}
f(X) = \bigcap \bigl\{F \in \cF \, : \, X \subset F \bigr\}, \ \ \ X
\subset \bL,
\end{equation}
where $\cF$ is a prescribed family of subsets of $\bL$.  The convex
hull itself and many well-known  generalizations of it are obtained
in this way, see for example \cite{bms} and \cite{v}. An immediate
consequence of definition \eqref{eqff} is the inclusion $X \subseteq
f(X)$.  From examples in Section~\ref{sdbp} and our results in
Section~\ref{scoh} it is clear that  the convex interval hull does
not always satisfy the inclusion $S \subseteq \co(S,\cJ_S)$.  As a
matter of fact, for every set $S$ there are families of intervals
$\cJ_S$ for which $S$ is not a subset of $\co(S,\cJ_S)$. In this
sense our operator differs from  operators described by
\eqref{eqff}.

Definitions similar to Definition~\ref{dco} appeared in \cite{gg}
and \cite{mo}. We recall the following three definitions from
\cite[p. 363]{mo}. First, for nonempty sets $\Lambda \subset \rr^m$
and $S \subset \bL$ denote by $\Lambda \cdot S \subset \bL$  the set
of all $\sum_{j=1}^m \lambda_j s_j$, where
$(\lambda_1,\ldots,\lambda_m) \in \Lambda$ and $s_j \in S, \,
j=1\ldots,m$.  Second, a set $S \subset \bL$ is called {\it
endo\,}-$\!\Lambda$ if $\Lambda \cdot S \subseteq S$. Third, with
$\cF$ being the family of all endo-$\!\Lambda$ sets, \eqref{eqff}
defines the $\Lambda$-hull operator.  A special case of
$\Lambda$-hull operator with $\Lambda = \bigl\{ (\xi,1-\xi) : \xi
\in \Delta \bigr\} \subset \nR^2$, where $\Delta$ is any non-empty
subset of $[0,1]$ containing at least one point interior to $[0,1]$,
was considered in \cite{gg}.  In \cite{gg} endo-$\Lambda$ sets are
called $\Delta$-convex sets. (We notice that Motzkin in \cite{mo}
does not refer to \cite{gg}.)

In this paragraph we point out the differences between the
definitions of $\Lambda \cdot S$ and $\co(S,\cJ_S)$.  To this end,
let $\Lambda$ be the intersection of $I_1\times \cdots \times I_m$
and the hyperplane $\sum_{j=1}^m \xi_j = 1$, where $I_j$ are
nonempty intervals in $\nR$, and let $S = \{x_1,\ldots,x_m\}$, where
$x_1,\ldots,x_m$ are distinct points in $\bL$.  Then, in general,
$\Lambda \cdot S$ contains more linear combinations than
$\co(S,J_S)$.  The first reason for this is that, with
$\bigl(\xi_1,\ldots\xi_m\bigr) \in \Lambda$,  $\sum_{j=1}^m \xi_j
s_j \in \Lambda \cdot S$ whenever $s_1,\ldots,s_m \in S$, while for
$\sum_{j=1}^m \xi_j x_j \in \co(S,\cJ_S)$ it is essential that
$x_1,\ldots,x_m$ are distinct points in $S$.  For example, with $s_1
= \cdots = s_m = s \in S$, the condition $\sum_{j=1}^m \xi_j = 1$
implies $S \subset \Lambda \cdot S$, while $S \subset \co(S,\cJ_S)$
is not true in general. The second reason is that in the definition
of $\co(S,\cJ_S)$ the point $x_j \in S$, for fixed $j \in
\{1,\ldots,m\}$, is scaled {\em only} by scalars in $I_j$, while
there is no such restriction in the definition of $\Lambda\cdot S$.
We also remark that the geometry of the sets $\Lambda \cdot S$ and
the properties of the operator $S \mapsto \Lambda\cdot S$ for a
fixed $\Lambda$ were not considered in \cite{gg} and \cite{mo}.

The article is organized as follows. In Section~\ref{sdbp} we give
several illustrative examples of convex interval hulls in $\nR^2$
and $\nR^3$ for sets $S$ with three, four and five points. A
justification for the adjective ``convex'' in the name of
$\co(S,\cJ_S)$  begins Section~\ref{sbp}. Furthermore, in this
section we characterize nonemptyness and boundedness of
$\co(S,\cJ_S)$.  In Section~\ref{spol} we prove that all bounded
convex interval hulls are polytopes.  We also give an upper bound
for the number of vertices of such polytopes.  As we have already
noticed, different families of intervals can result in the same
convex interval hulls. In Section~\ref{smi} we study minimality
conditions for a family of intervals, where  minimality is
understood in such a way that any further shrinking of intervals
results in a smaller convex interval hull. In Section~\ref{scoh} we
prove that a family of bounded convex interval hulls of a fixed
finite set $S$ is invariant under homotheties. As a special case of
this result we obtain that for each homothet $K$ of $\conv S$ there
exists a family of intervals $\cJ_S$ such that $K = \co(S,\cJ_S)$.
We use this result to give a detailed description of bounded convex
interval hulls of finite affinely independent subsets of a linear
space.

In this paragraph we introduce the notation. By $\nR$ we denote the
real numbers. The symbol $\bL$ denotes a real linear space and
$\|\cdot\|$ is a norm in this space.  A specific linear space that
we will encounter is $\nR^m$, where $m$ is a positive integer. The
linear operations from $\bL$ are extended to subsets of $\bL$ in the
following standard way. For subsets $K$ and $M$ of $\bL$ and
$\alpha, \beta \in \nR$ we put
 \[
\alpha K+ \beta M = \bigl\{\alpha x+ \beta y \; :\; x \; \in K,\;
y\in M \bigr\}.
 \]
For a mapping $T: \bL \to \bL$, $T(K)$ denotes the set of all $Tx,
\, x\in K$.

\section{Examples} \label{sdbp}

In this section we present several examples of convex interval
hulls.  All examples here are bounded sets, since our main interest
in this article are convex interval hulls which are bounded sets. We
will consider unbounded convex interval hulls elsewhere.  For
completeness we start with the standard example.

\begin{example} \label{ex1}
Let $S = \{x_1,\ldots,x_m\}$ be a finite set of points in a linear
space $\bL$ and let $\cJ_S = \{ I_1, \ldots, I_m\}$ with $I_j =
[0,t_j]$, where $t_j \geq 1,  \, j = 1,\ldots,m$. Then
 \[
\co(S,\cJ_S) = \conv S.
 \]
\end{example}

The examples below are calculated and plotted using {\em
Mathematica}.  In each example the points of the set $S$ are listed
starting from the lowest point that is furthermost to the left. Then
we proceed counterclockwise, finishing with the point inside. In
each figure the points in $S$ are marked with black dots ($\bullet$)
and the polygon $\co\bigl(S,\cJ_S\bigr)$ is shaded gray with its
edges slightly darker.

 \setlength{\abovecaptionskip}{0pt}%


\begin{figure}
\begin{minipage}[b]{.495\linewidth}
\resizebox{\linewidth}{!}{%
  \includegraphics{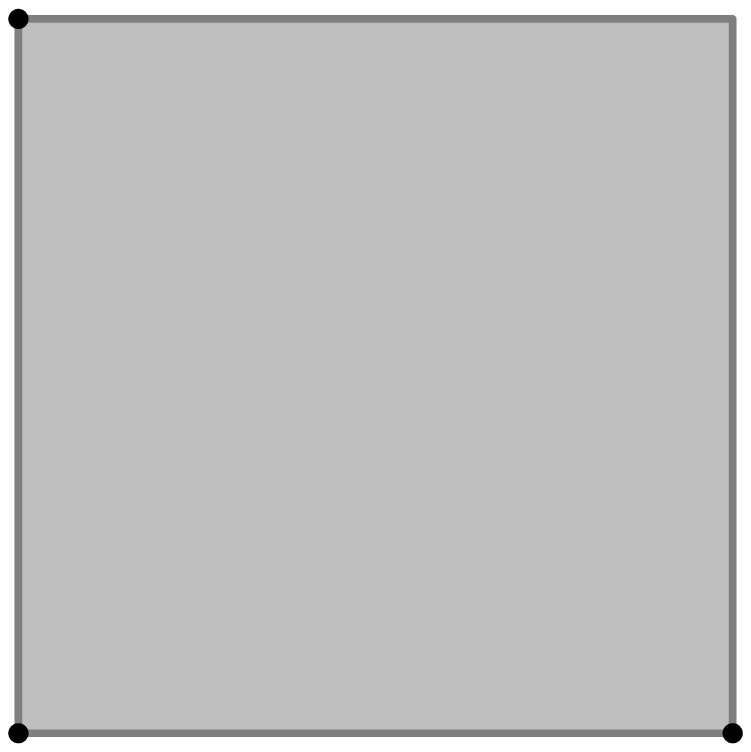}}
    \caption{A square}
\label{f1}
\end{minipage}
\hfill
\begin{minipage}[b]{.495\linewidth}
\resizebox{\linewidth}{!}{%
  \includegraphics{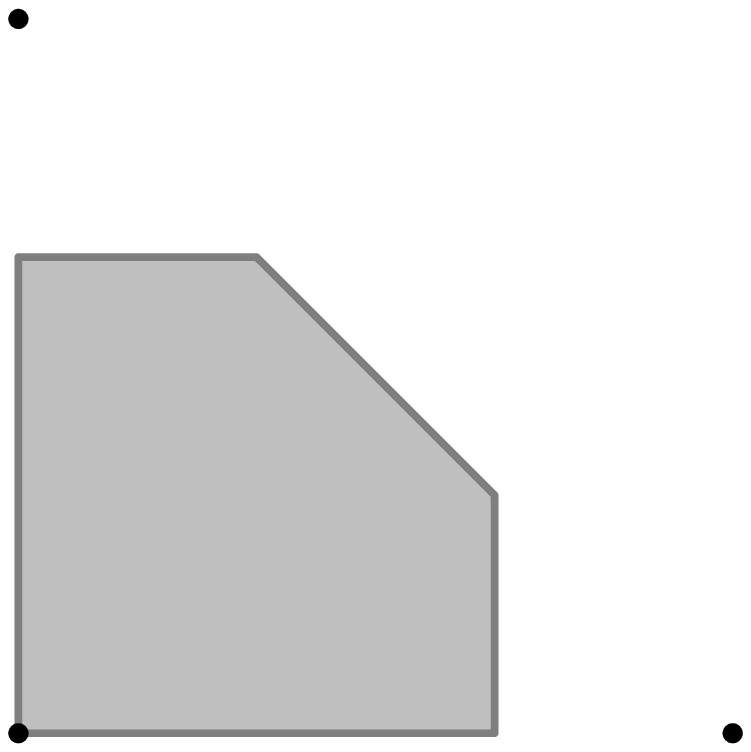}}
    \caption{A pentagon}
 \label{f2}
\end{minipage}


\begin{minipage}[b]{.495\linewidth}
\resizebox{\linewidth}{!}{%
  \includegraphics{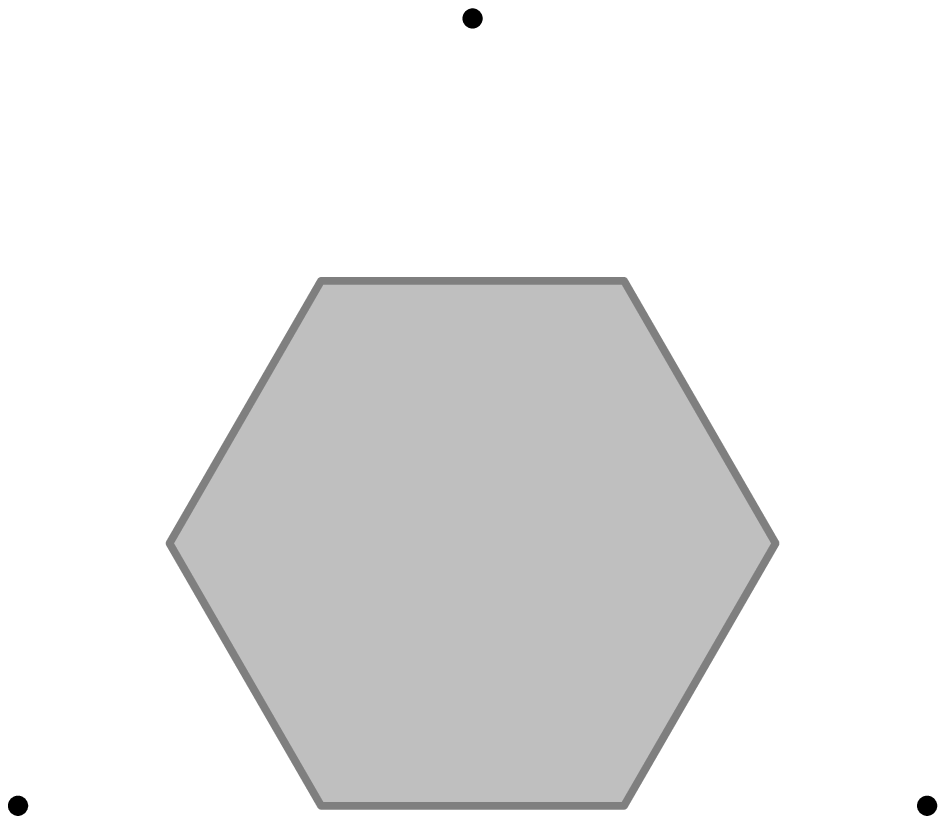}}
    \caption{A regular hexagon}
 \label{f3}
\end{minipage}
\hfill
\begin{minipage}[b]{.495\linewidth}
\resizebox{\linewidth}{!}{%
  \includegraphics{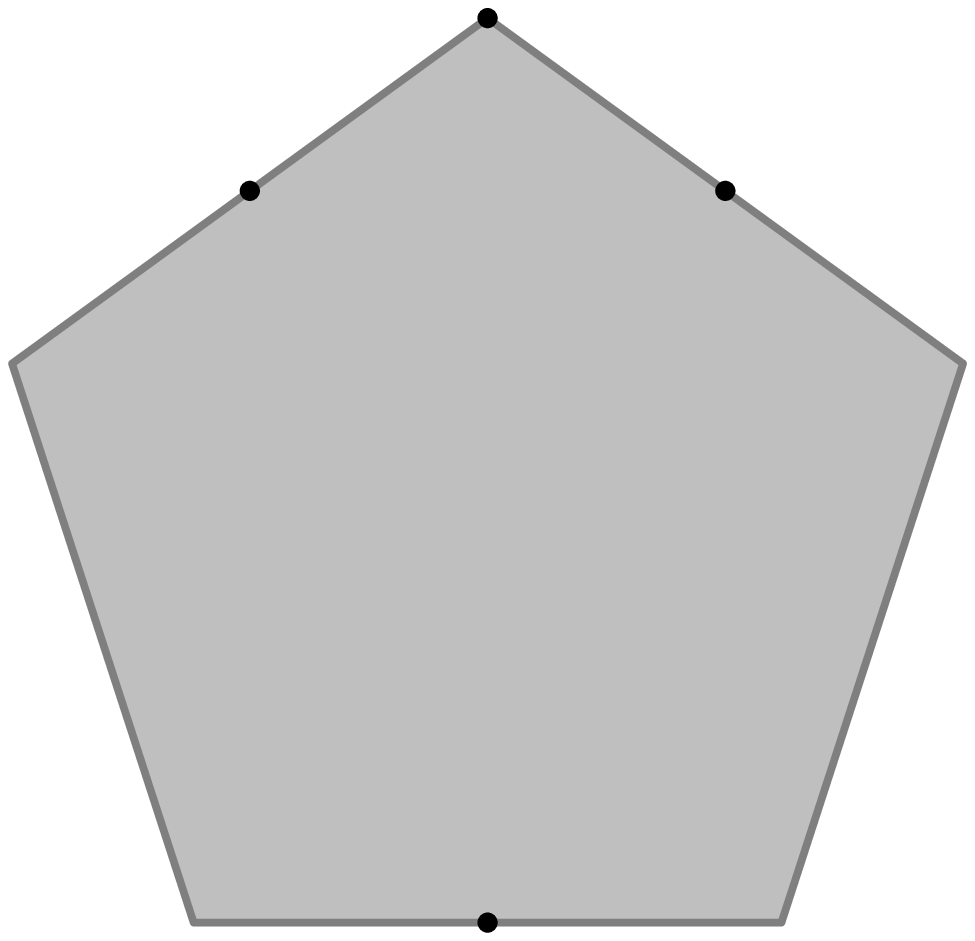}}
    \caption{A regular pentagon}
 \label{f4}
\end{minipage}


\begin{minipage}[b]{.495\linewidth}
\resizebox{\linewidth}{!}{%
  \includegraphics{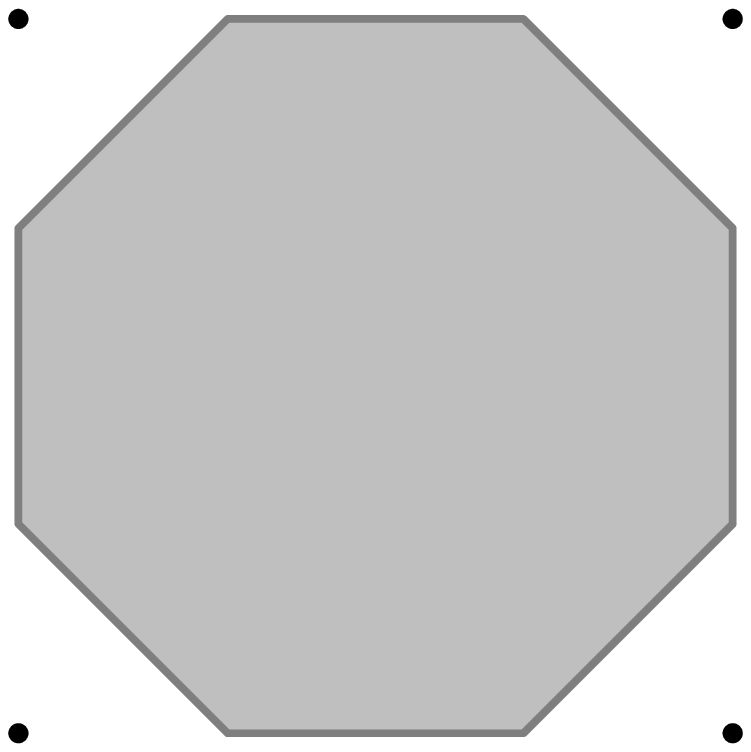}}
    \caption{A regular octagon}
 \label{f5}
\end{minipage}
\hfill
\begin{minipage}[b]{.495\linewidth}
\resizebox{\linewidth}{!}{%
  \includegraphics{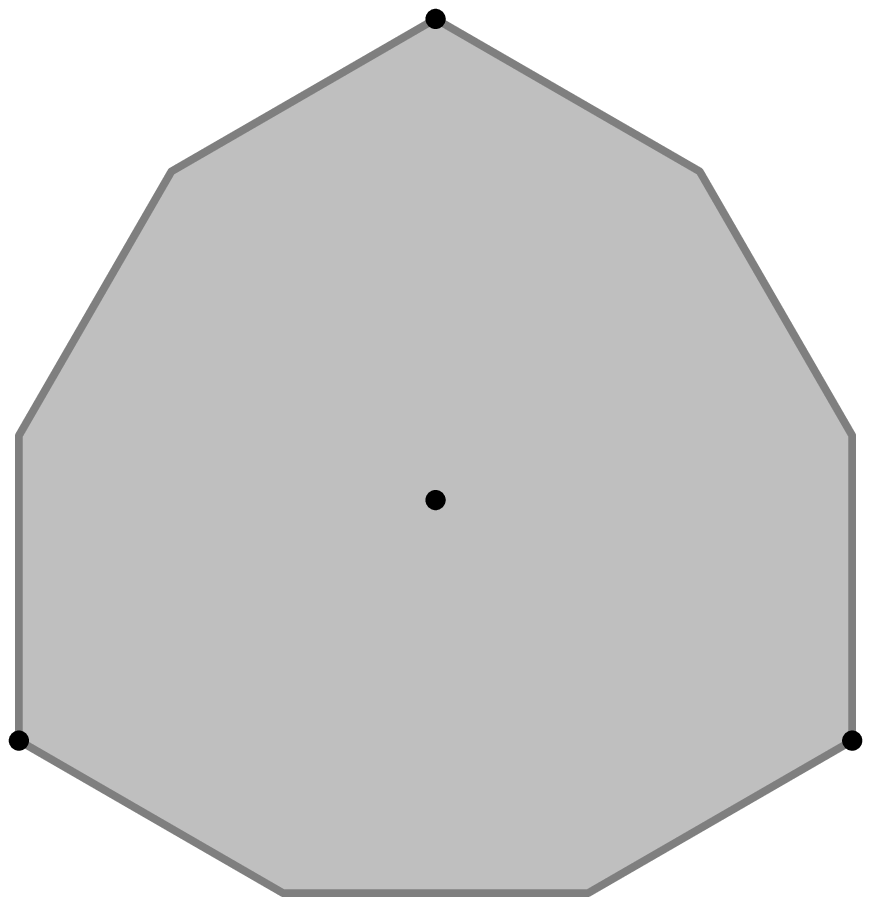}}
    \caption{A nonagon}
 \label{f6}
\end{minipage}
\end{figure}

\begin{example} \label{ex22}
In Figures~\ref{f1} and \ref{f2}, we use
$S=\bigl\{(0,0),(1,0),(0,1)\bigr\}$.  In Fig.~\ref{f1} we use
$\cJ_S=\bigl\{[-1,1],[0,1],[0,1]\bigr\}$ to get a square and in
Fig.~\ref{f2} we use $\cJ_S=\bigl\{[0,1],[0,2/3],[0,2/3]\bigr\}$ to
get an irregular pentagon.
\end{example}

\begin{example} \label{ex23}
In Fig.~\ref{f3} we use
\[
S=\bigl\{(-1,0),(0,1),(0,\sqrt{3})\bigr\} \ \ \ \text{and} \ \ \
\cJ_S = \bigl\{ [0,2/3], [0,2/3], [0,2/3] \bigr\},
\]
to get a regular hexagon.
\end{example}

\begin{example}
In Fig.~\ref{f4} we use
\begin{align*}
 S  & = \left\{\!
   \Bigl( 0,-5 - 2 {\sqrt{5}}\, \Bigr),
  \Bigl( \sqrt{10+2\sqrt{5}},\sqrt{5} \, \Bigr) ,
  \bigl(0,5 \bigr) ,
  \Bigl( - \sqrt{10+2\sqrt{5}},\sqrt{5}\, \Bigr)
\right\}, \\
 \cJ_S & = \bigl\{[0,3 - \sqrt{5}], [0,2], [-1,1], [0,2] \bigr\},
\end{align*}
to get a regular pentagon.
\end{example}

\begin{example}
In Fig.~\ref{f5} we use $S = \bigl\{(0,0),(1,0),(1,1),(0,1)\bigr\}$
and $\cJ_S$ that consists of four copies of the interval
$\bigl[0,\sqrt{2}/2\bigr]$ to get a regular octagon.
\end{example}

\begin{example} \label{ex4pts}
In Fig.~\ref{f6} we use
\begin{align*}
S & = \bigl\{\bigl(-1,0\bigr), \bigl(0,1\bigr),
\bigl(0,\sqrt{3}\,\bigr),
\bigl( 0,1/\sqrt{3}\, \bigr) \bigr\}, \\
\cJ_S & = \bigl\{ [0,1], [0,1], [0,1],\bigl[(\sqrt{3}-3)/2,1\bigr]
\bigr\},
\end{align*}
to get an irregular nonagon with all equal sides.
\end{example}

\begin{figure}
\begin{minipage}[b]{.495\linewidth}
\resizebox{\linewidth}{!}{%
  \includegraphics{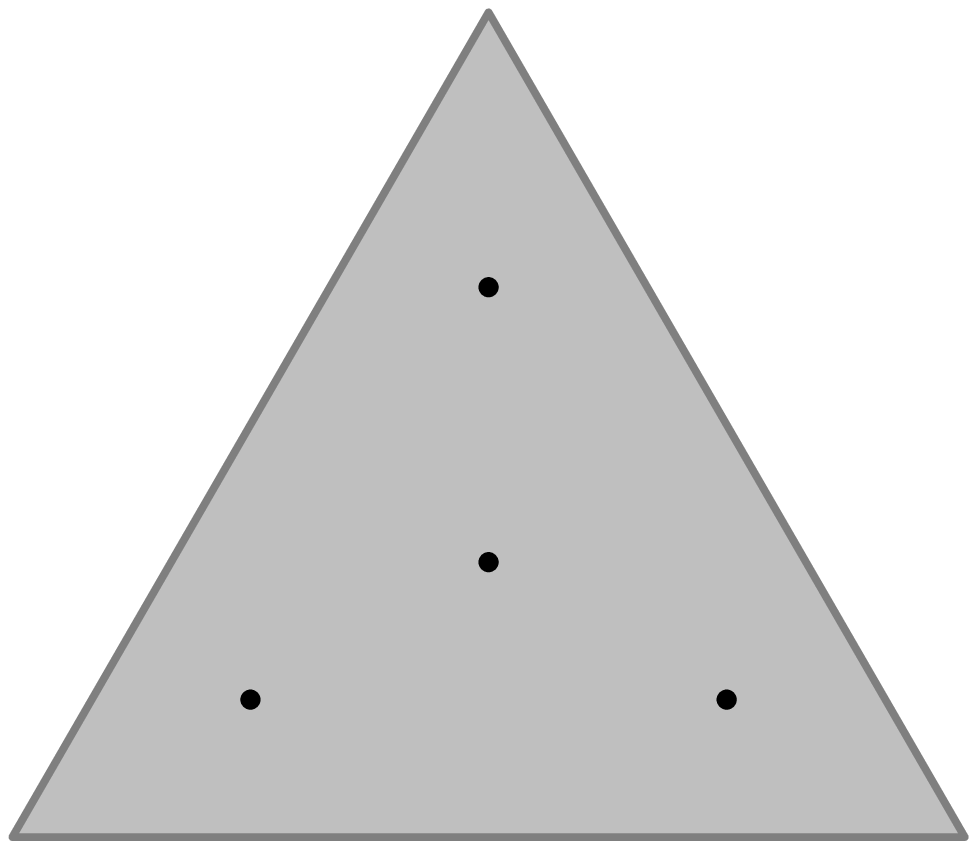}}
    \caption{Step 1}
 \label{f7}
\end{minipage}
\hfill
\begin{minipage}[b]{.495\linewidth}
\resizebox{\linewidth}{!}{%
  \includegraphics{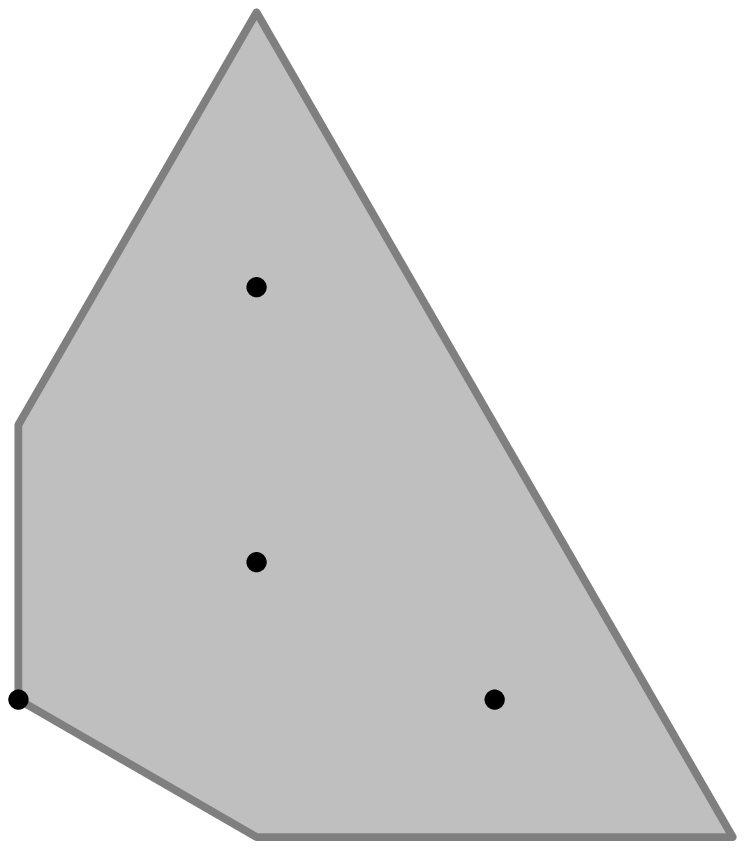}}
    \caption{Step 2}
 \label{f8}
\end{minipage}


\begin{minipage}[b]{.495\linewidth}
\resizebox{\linewidth}{!}{%
  \includegraphics{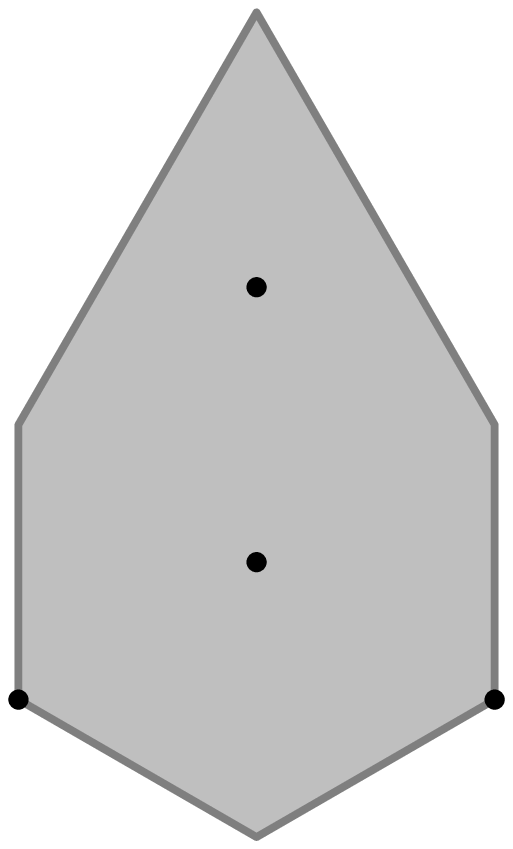}}
    \caption{Step 3}
 \label{f9}
\end{minipage}
\hfill
\begin{minipage}[b]{.495\linewidth}
\resizebox{\linewidth}{!}{%
  \includegraphics{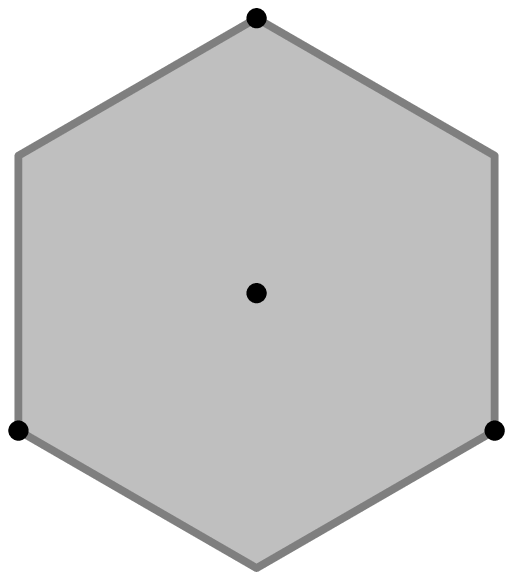}}
    \caption{Step 4}
 \label{f10}
\end{minipage}


\begin{minipage}[b]{.495\linewidth}
\resizebox{\linewidth}{!}{%
  \includegraphics{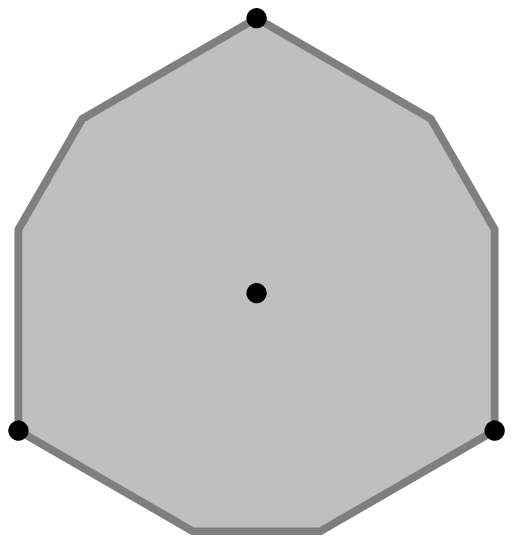}}
    \caption{Step 5}
 \label{f11}
\end{minipage}
\hfill
\begin{minipage}[b]{.495\linewidth}
\resizebox{\linewidth}{!}{%
  \includegraphics{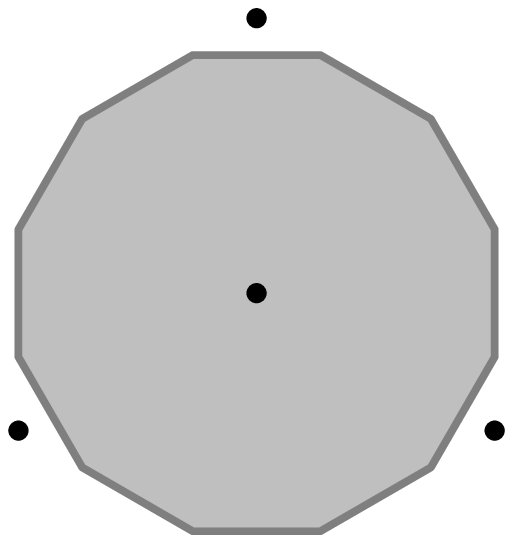}}
    \caption{A regular dodecagon}
 \label{f12}
\end{minipage}
\end{figure}

\begin{example} \label{ex3}
In Figures~\ref{f7} through~\ref{f12} we show six different convex
interval hulls corresponding to the same set $S$ that is used in
Example~\ref{ex4pts}. We start with an equilateral triangle in
Fig.~\ref{f7} and proceed by changing one interval at each step to
finish with a regular dodecagon in Fig.~\ref{f12}. We use the
following families of intervals:
\begin{alignat*}{2}
&\text{Fig.~\ref{f7}} \qquad \qquad
 & \cJ_S & = \bigl\{ [0,2], [0,2], [0,2], [-1,0] \bigr\}, \\
&\text{Fig.~\ref{f8}}
 & \cJ_S & = \bigl\{ [0,1], [0,2], [0,2], [-1,0] \bigr\}, \\
&\text{Fig.~\ref{f9}}
 & \cJ_S & = \bigl\{ [0,1], [0,1], [0,2], [-1,0] \bigr\}, \\
&\text{Fig.~\ref{f10}}
 & \cJ_S & = \bigl\{ [0,1], [0,1], [0,1], [-1,0] \bigr\}, \\
&\text{Fig.~\ref{f11}}
 & \cJ_S & = \bigl\{ [0,1], [0,1], [0,1], [1-\sqrt{3},0] \bigr\}, \\
&\text{Fig.~\ref{f12}}
 & \cJ_S & = \bigl\{ [0,1], [0,1], [0,1],
     [1-\sqrt{3},-2+\sqrt{3} \,] \bigr\}.
\end{alignat*}
\end{example}



\begin{figure}
\begin{minipage}[b]{.495\linewidth}
\resizebox{\linewidth}{!}{%
  \includegraphics{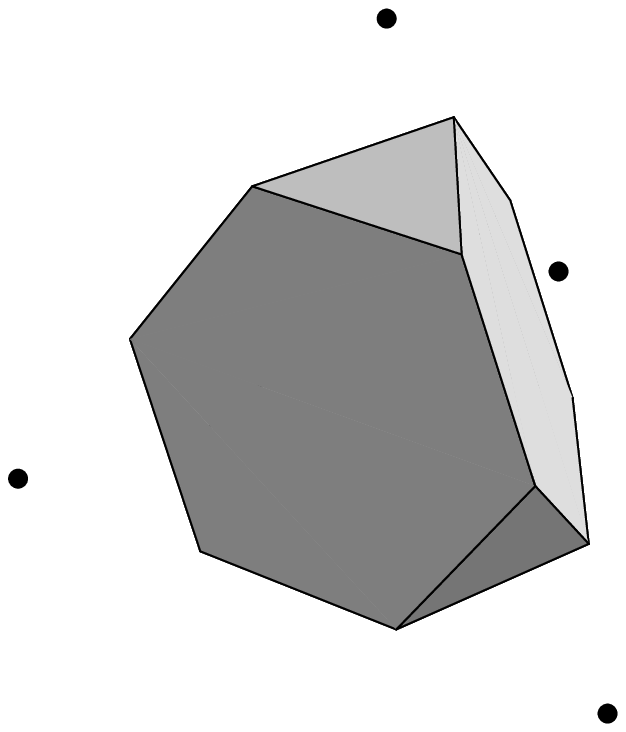}}
    \caption{A~truncated~tetrahedron}
 \label{f13}
\end{minipage}
 \hfill
\begin{minipage}[b]{.495\linewidth}
\resizebox{\linewidth}{!}{%
  \includegraphics{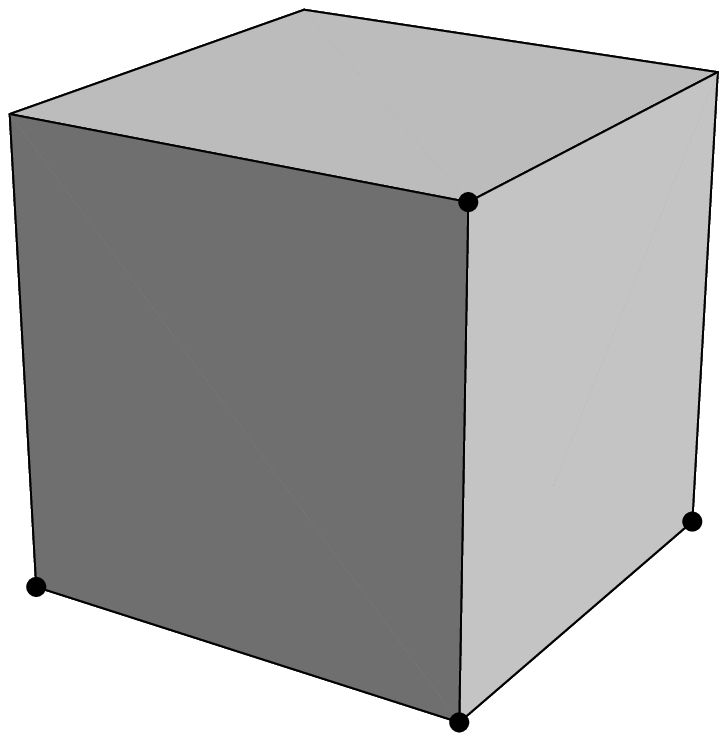}}
    \caption{A cube}
 \label{f14}
\end{minipage}

\vspace{5pt}

\begin{minipage}[b]{.495\linewidth}
\resizebox{\linewidth}{!}{%
  \includegraphics{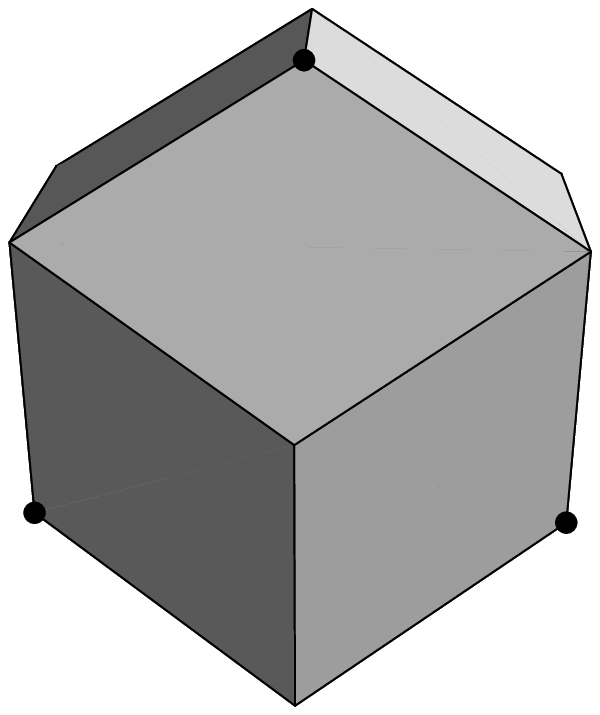}}
    \caption{A rhombic dodecahedron}
 \label{f15}
\end{minipage}
\hfill
\begin{minipage}[b]{.495\linewidth}
\resizebox{\linewidth}{!}{%
  \includegraphics{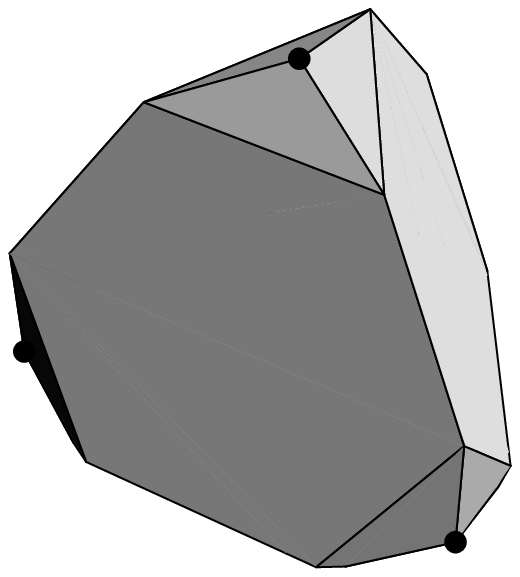}}
    \caption{Example~\ref{exttps}}
 \label{f16}
\end{minipage}

\end{figure}
%
%

\begin{example} \label{ex28}
In Fig. \ref{f13} we use
\begin{align*}
S & = \Bigl\{ \bigl( -1,0,0 \bigr) ,\bigl( 1,0,0 \bigr) ,
  \bigl( 0,{\sqrt{3}},0 \bigr) ,
  \bigl( 0,1/\sqrt{3},2\,{\sqrt{{2}/{3}}} \,\bigr) \Bigr\}
\end{align*}
and $\cJ_S$ that consists of four copies of $\bigl[0,2/3\bigr]$ to
get a truncated tetrahedron. Notice that the points of $S$ are
vertices of a tetrahedron.
\end{example}

\begin{example} \label{ex29}
In Fig. \ref{f14} we use
\begin{align*}
S & = \Bigl\{ \bigl( 0,0,0 \bigr) ,\bigl( 1,0,0 \bigr) ,
  \bigl( 0,1,0 \bigr) ,
  \bigl( 0,0,1 \bigr) \Bigr\} \\
\cJ_S & = \Bigl\{ \bigl[-2,1\bigr],\bigl[0,1\bigr],
  \bigl[0,1\bigr] ,
  \bigl[0,1\bigr] \Bigr\}
\end{align*}
to get a cube.
\end{example}

\begin{example} \label{exrhdo}
In Fig. \ref{f15} we use
\begin{align*}
S & = \left\{\! \Bigl( -1,0,0 \Bigr) ,\Bigl( 1,0,0 \Bigr) ,
  \Bigl( 0,\sqrt{3},0 \Bigr) ,
  \left( 0,\frac{1}{\sqrt{3}},2\sqrt{\frac{2}{3}}\right),
  \left( 0,\frac{1}{\sqrt{3}},\frac{1}{\sqrt{6}}\right)  \!\right\} \\
\cJ_S & = \Bigl\{ \bigl[0,1\bigr],\bigl[0,1\bigr],
  \bigl[0,1\bigr] ,
  \bigl[0,1\bigr],\bigl[-2,0\bigr] \Bigr\}
\end{align*}
to get a rhombic dodecahedron. The first four points of $S$ are
vertices of a tetrahedron and the fifth point is its orthocenter.
\end{example}

\begin{example} \label{exttps}
In Fig. \ref{f16} we use the same $S$ as in Example~\ref{exrhdo} and
\begin{align*}
\cJ_S & = \Bigl\{ \bigl[0,1\bigr],\bigl[0,1\bigr],
  \bigl[0,1\bigr] ,
  \bigl[0,1\bigr],\bigl[-1/2,0\bigr] \Bigr\}.
\end{align*}
\end{example}

\section{Basic properties of convex interval hulls} \label{sbp}

In this section most proofs are omitted since they are, though
sometimes lengthy, straightforward consequences of the definitions.
The proofs that are included indicate how to construct the omitted
proofs.

\begin{proposition} \label{pco}
Let $S = \{x_1,\ldots,x_m\}$ be a finite set of points in a linear
space $\bL$ and let $\cJ_S = \{ I_1, \ldots, I_m\},\, I_j \subseteq
\nR, \, j = 1,\ldots,m$, be a family of nonempty intervals. Then the
set $\co(S,\cJ_S)$ is convex.
\end{proposition}

\begin{proposition} \label{ppco1}
Let $S = \{x_1,\ldots,x_m\}$ be a finite set of points in a linear
space $\bL$ and let $\cJ_S = \{ I_1, \ldots, I_m\},\, I_j \subseteq
\nR, \, j = 1,\ldots,m$, be a family of nonempty intervals.  If
$\cJ_S^{\prime} = \{ I_1^{\prime}, \ldots, I_m^{\prime}\}$ is a
family of nonempty intervals such that $I_j^{\prime} \subseteq I_j,
\, j=1,\ldots,m$, then
 \[
\co(S,{\cJ}_S^{\prime})\subseteq \co(S,{\cJ}_S).
 \]
\end{proposition}

\begin{proposition} \label{ppco2}
Let $S = \{x_1,\ldots,x_m\}$ be a finite set of points in a linear
space $\bL$ and let $\cJ_S = \{ I_1, \ldots, I_m\},\, I_j \subseteq
\nR, \, j = 1,\ldots,m$, be a family of nonempty intervals. Put $a_j
= \inf I_j$ and $b_j = \sup I_j, \, \jom$, allowing for the infinite
values. Let $\alpha = \su a_j$ and $\beta = \su b_j$. Then
$\co(S,\cJ_S) \neq \emptyset$ if and only if the following three
conditions are satisfied:
\begin{enumerate}[\rm (a) \ ]
\item
$\alpha \leq 1 \leq \beta$;
\item
if $\alpha = 1$, then $a_j \in I_j, \, \jom$; and
\item
if $\beta = 1$, then $b_j \in I_j, \, \jom$.
\end{enumerate}
\end{proposition}

\begin{proof}
Assume $\co(S,\cJ_S) \neq \emptyset$.  Then there exist $\xi_j \in
I_j, \, \jom$, such that $\su \xi_j =1$. Since $a_j \leq \xi_j \leq
b_j, \, \jom$, it follows that $\alpha \leq 1 \leq \beta$. This
proves (a). If $\alpha = 1$, then $a_j = \xi_j$, and thus $a_j \in
I_j, \ \jom$. This proves (b) and (c) is proved similarly.

To prove the converse assume that (a), (b) and (c) hold.  If $\alpha
= \beta = 1$, then each of the intervals is in fact a point and
$\co(S,\cJ_S)$ consists of a single point. Now assume $\alpha <
\beta$. It follows from Proposition~\ref{ppco1} that without loss of
generality we can in addition assume that $\alpha$ and $\beta$ are
finite. Set
 \[
\xi_j = \dfrac{\beta - 1}{\beta - \alpha}\, a_j +
 \dfrac{1-\alpha}{\beta - \alpha}\, b_j, \quad \jom.
 \]
It easily follows that $\xi_j \in I_j,\, \jom$, and $\su \xi_j = 1$.
Therefore
 $
x = \su \xi_j x_j \in \co(S,\cJ_S).
 $
Thus $\co(S,\cJ_S) \neq \emptyset$ and the proposition is proved.
\end{proof}

\begin{theorem} \label{tpco}
Let $S = \{x_1,\ldots,x_m\}$ be a finite set of distinct points in a
linear space $\bL$ and let $\cJ_S = \{ I_1, \ldots, I_m\},\, I_j
\subseteq \nR, \, j = 1,\ldots,m$, be a family of nonempty
intervals.

The set $\co(S,\cJ_S)$ is bounded if and only if at least one of the
conditions below is satisfied.
\begin{enumerate}[{\rm (i)}]
\item \label{tpcodiii}
All the intervals in $\cJ_S$ are bounded below.
\item \label{tpcodii}
All the intervals in $\cJ_S$ are bounded above.
\item \label{tpcodi}
At most one interval in $\cJ_S$ is unbounded.
\end{enumerate}

If any of the conditions {\rm (i)-(iii)} is satisfied, then there
exists a family of bounded intervals $\cJ_S^{\prime}$ such that
$\co(S,\cJ_S) = \co(S,\cJ_S^{\prime})$.
\end{theorem}

\begin{proof}
We first show that if all the intervals in $\cJ_S =
\{I_1,\ldots,I_m\}$ are bounded, then $\co(S,\cJ_S)$ is bounded.
Indeed, if for some $b > 0$ we have $I_j \subseteq [-b,b]$ for all
$\jom$, then for each $x \in \co(S,\cJ_S)$ we have
 \[
\|x\| = \biggl\| \, \su \xi_j\, x_j \, \biggr\| \leq m\, b\,
\max\bigl\{ \|x_1\|,\ldots, \|x_m\|\bigr\}.
 \]

Next we assume (\ref{tpcodiii}).  Let $a \in \nR$ be such that $I_j
\subseteq (a,+\infty)$ for all $\jom$.  Since the empty set is
bounded, we assume $\co(S,\cJ_S) \neq \emptyset$.  By
Proposition~\ref{ppco2}, we have $m a < 1$. Let $x \in
\co(S,\cJ_S)$. Then there exist $\xi_j \in I_j,\, \jom,$ such that
$\su \xi_j =1$ and $x =\su \xi_j\, x_j$. Let $k \in \{1,\ldots,m\}$
be arbitrary. Then
 \[
\xi_k = 1 - \sum\limits_{j\not=k}^{m} \xi_j
 < 1- (m-1)a = 1-m a + a.
 \]
Therefore
 \[
a < \xi_k < 1 -  m a + a.
 \]
Hence, with $b_k^{\prime} = \min \{1 -  m a + a,b_k\}$,
$I_k^{\prime} = [a_k,b_k^{\prime}]$ and $\cJ_S^{\prime} =
\{I_1^{\prime}, \ldots, I_m^{\prime} \}$, we have proved that
$\co(S,\cJ_S) \subseteq \co(S,\cJ_S^{\prime})$.  By
Proposition~\ref{ppco1} the converse inclusion is also true.
Consequently $\co(S,\cJ_S) = \co(S,\cJ_S^{\prime})$. Since each
interval in $\cJ_S^{\prime}$ is bounded, the set
$\co(S,\cJ_S^{\prime})$ is bounded. Thus $\co(S,\cJ_S)$ is bounded,
as well.

Similarly, (\ref{tpcodii}) implies that $\co(S,\cJ_S)$ is bounded.

Assume (\ref{tpcodi}). We can also assume that (\ref{tpcodiii}) and
(\ref{tpcodii}) are not true.  Then exactly one of the intervals in
$\cJ_S$ is unbounded and it equals $\nR$.  Assume that $I_1 = \nR$
and $I_j = [a_j,b_j], a_j \leq b_j, a_j, b_j \in \nR, j=2,\ldots,m$.
Put $\alpha = \sum_{j=2}^m a_j$ and $\beta = \sum_{j=2}^m b_j$.
Clearly $\alpha \leq \beta$.  Let $v \in \co(S,\cJ_S)$ be such that
 \[
v = \su \xi_j x_j, \ \ \ \su \xi_j = 1, \ \ \ \xi_j \in I_j, \ \ \
\jom.
 \]
Then
\[
1-\beta \leq \xi_1 = 1 - \xi_2 - \cdots - \xi_m \leq 1 - \alpha.
\]
Consequently $v \in \co(S,\cJ^{\prime}_{S})$, where
$\cJ^{\prime}_{S} = \bigl\{[1-\beta,1-\alpha],I_2,\ldots,I_m
\bigr\}$. Therefore $\co(S,\cJ_{S}) \subseteq
\co(S,\cJ^{\prime}_{S})$.  The converse inclusion holds by
Proposition~\ref{ppco1}. Therefore $\co(S,\cJ_{S}) =
\co(S,\cJ^{\prime}_{S})$. Since each interval in $\cJ^{\prime}_{S}$
is bounded, the set $\co(S,\cJ^{\prime}_{S})$ is bounded and so is
$\co(S,\cJ_{S})$.  This completes the proof of "if'' part of the
theorem.

Next we prove the contrapositive of the ``only if'' part of the
theorem.  Assume that (\ref{tpcodiii}), (\ref{tpcodii}) and
(\ref{tpcodi}) are all false.  This is equivalent to the fact that
the family $\cJ_S$ contains at least two unbounded intervals, say
$I_1$ and $I_2$, such that $I_1$ is not bounded from below and $I_2$
is not bounded from above.  Let $v \in \co(S,\cJ_S)$ be such that
 \[
v = \su c_j x_j, \ \ \ \su c_j = 1, \ \ \ c_j \in I_j, \ \ \ \jom.
 \]
Then
\[
(-\infty,c_1] \subseteq I_1  \ \ \ \text{and} \ \ \ [c_2, +\infty)
\subseteq I_2.
 \]
Consequently, for all $t \geq 0$,
 \[
(c_1-t) x_1 + (c_2 + t) x_2 + c_3 x_3 + \cdots + c_m x_m = v +
t(x_2-x_1) \in \co(S,\cJ_S).
 \]
Clearly
 \[
\|v + t(x_2-x_1)\| \geq t\|x_2-x_1\| - \|v\|, \ \ \ t \geq 0.
 \]
Since by assumption $x_1 \neq x_2$, the last inequality implies that
$\co(S,\cJ_S)$ is unbounded. The theorem is proved.
\end{proof}

\begin{proposition} \label{p1e}
Let $T:\bL \to \bK$ be an affine transformation between linear
spaces $\bL$ and $\bK$. Let $S = \bigl\{x_1,\ldots,x_m \bigr\}$ be a
finite subset of $\bL$ and let $\cJ_S =
\bigl\{I_1,\ldots,I_m\bigr\}$ be a corresponding set of intervals
for which $\co(S,\cJ_S)$ is bounded. Let $Q = T(S) = \{y_1,\ldots,
y_k\}$ be the set with $k$ elements, where $k \leq m$. Set
\begin{equation*}
\cJ_Q = \bigl\{ I_1^{\prime}, \ldots, I_k^{\prime} \bigr\} \quad
\text{where} \quad I_j^{\prime} := \sum \bigl\{ I_i \, : \, T x_i =
y_j \bigr\}.
\end{equation*}
Then
\begin{equation*} 
\co\bigl(Q, \cJ_Q\bigr) = T \bigl(\co(S,\cJ_S)\bigr).
\end{equation*}
Each vertex of $T \bigl(\co(S,\cJ_S)\bigr)$ is an image of a vertex
of $\co(S,\cJ_S)$.
\end{proposition}

\begin{proposition}
Let  $S$ be a finite subset of \ $\bL$ and let $\,T: {\bL} \to \bL$
be an affine bijection such that $T(S) = S$.  Assume that ${\cJ}_S$
has the property that $I_j = I_k$ whenever $x_j = T x_k$. Then
$T\bigl(\co(S,{\cJ}_S)\bigr)  = \co(S,{\cJ}_S)$.
\end{proposition}

\begin{corollary}
Let  $S \subset \bL$ be a finite set centrally symmetric with
respect to $u\in\bL$. Assume that $\cJ_S$ has the property that
$I_i=I_j$ for $x_i$ and its symmetric image $x_j$. Then
$\co(S,\cJ_S)$ is also centrally symmetric with respect to $u$.
\end{corollary}

\section{Convex interval hulls and polytopes} \label{spol}

Example~\ref{ex3} shows that a convex interval hull of four points
can have twelve vertices. In the next theorem we give an upper bound
for the number of vertices of a convex interval hull for a finite
set with $m$ points. For a real number $t$, $\lfloor t \rfloor$
denotes the greatest integer that does not exceed $t$.

\begin{theorem}\label{tmax}
Let $S = \{x_1,\ldots,x_m\}$ be a subset of $\,\bL$ and let $\cJ_S$
be a family of closed intervals such that $\co(S,\cJ_S)$ is bounded.
Then $\co(S,\cJ_S)$ is the  convex hull of at most
\begin{equation*}
n \binom{m}{n}  \quad \text{points, \ where} \quad n = \Bigg\lfloor
\dfrac{m}{2} \Bigg\rfloor + 1
\end{equation*}
and this bound is best possible.
\end{theorem}
\begin{proof}
It follows from Proposition~\ref{ppco2} that there is no loss of
generality if we assume that all the intervals in
 $\cJ_S$ are bounded. Set $I_j =[a_j,b_j], \, a_j
< b_j, \, j=1,\ldots,m$, and, as before, set $\alpha = \sum_{j=1}^m
a_j$.  Clearly, for $\alpha = 1$ the set $\co(S,\cJ_S)$ is a
singleton and the theorem is trivial in this case. Therefore,  we
can assume that $\alpha < 1$.

In order  to prove the theorem we will show that $\co(S,\cJ_S)$ is
an image under an affine transformation of a polytope in $\rr^m$
having not more than $n \binom{m}{n}$ vertices. To this end take
 the unit vectors $e_1, \ldots,e_m$  in $\nR^m$ and put
\begin{equation*}
Q = \left\{\dfrac{1-\alpha}{b_1 - a_1}\, e_1,
\ldots,\dfrac{1-\alpha}{b_m - a_m}\, e_m \right\}
\end{equation*}
assigning  to  $Q$ the family of intervals
 \[
\cJ_Q = \bigl\{I_1^{\prime},\ldots,I_m^{\prime} \bigr\}, \quad
\text{where} \quad I_j^{\prime} = \left[0,\dfrac{b_j -
a_j}{1-\alpha} \right], \,\, j=1,\dots,m.
 \]
For such $Q$ and $\cJ_Q$ we have
 \[
 \co\bigl(Q, \cJ_Q\bigr) = \left\{\sum\limits_{j=1}^m \eta_j\,
 \dfrac{1-\alpha}{b_j - a_j}\, e_j \, : \,
 \eta_j \in I_j^{\prime}, \, \,
 \sum\limits_{j=1}^m \eta_j = 1 \right\}.
 \]
This, after a straightforward substitution $\eta_j
\,\dfrac{1-\alpha}{b_j - a_j} = \mu_j$,  gives
\begin{align*}
\co\bigl(Q, \cJ_Q\bigr)& = \left\{\sum\limits_{j=1}^m \mu_j e_j \,
\, : \, 0\leq \mu_j \leq 1 \, , \,\,\,
\sum\limits_{j=1}^m(b_j -a_j)\, \mu_j = 1-\alpha \right\}\\
& = \left\{ (\mu_1,\ldots, \mu_m) \in C \, : \,\,\, \sum_{j=1}^m
(b_j - a_j)\, \mu_j = 1-\alpha \right\},
\end{align*}
where $C$ is the unit hypercube in $\nR^m$. From the above it is
immediately seen that $\co\bigl(Q, \cJ_Q\bigr) = D$ is the
intersection of the hypercube $C$ and the hyperplane
\begin{equation*} 
\Pi := \left\{ (\mu_1,\ldots, \mu_m) \in \nR^m \, : \, \sum_{j=1}^m
(b_j - a_j)\,  \mu_j = 1-\alpha \right\}.
\end{equation*}
Therefore, if  $v$ is a vertex of $D$, then $v$ must be on an edge of
$C$. Clearly, if $e$ is an edge of $C$
intersecting $\Pi$ but not lying on $\Pi$, then there can be only one vertex of $D$ on $e$.
By \cite[Theorem 5]{az} a hyperplane can intersect at most $n
\binom{m}{n}$ edges of $C$.
 Consequently, from the two observations it follows
that $D$ has at most $n \binom{m}{n}$ vertices.

Now consider the affine transformation $T: \nR^m \to \bL$ defined as
follows
\begin{equation*}
T : \sum\limits_{j=1}^m \zeta_j e_j \mapsto   \sum\limits_{j=1}^m
a_j x_j  + \sum\limits_{j=1}^m (b_j - a_j)\, \zeta_j x_j, \quad
(\zeta_1,\ldots, \zeta_m) \in \nR^m.
\end{equation*}
Calculating $T(D)$ we obtain
\begin{align*}
T(D)& = T \left( \Biggl\{\sum\limits_{j=1}^m \eta_j\,
\dfrac{1-\alpha}{b_j - a_j}\, e_j \, : \,\ \eta_j \in I_j^{\prime},
\,\, \sum\limits_{j=1}^m \eta_j = 1
\Biggr\} \right)  \\
 & = \Biggl\{ \sum\limits_{j=1}^m a_j x_j + \sum\limits_{j=1}^m
 (1-\alpha) \, \eta_j \, x_j \, : \,\ \eta_j \in I_j^{\prime}, \,\,
 \sum\limits_{j=1}^m \eta_j = 1  \Biggr\} \\
  & = \Biggl\{\sum\limits_{j=1}^m
  \bigl( a_j + (1-\alpha) \, \eta_j  \bigr)
 \, x_j \, : \,\ \eta_j \in I_j^{\prime}, \,\,
 \sum\limits_{j=1}^m \eta_j = 1 \Biggr\}.
\end{align*}
Substituting $\xi_j =  a_j + (1-\alpha) \, \eta_j, \, j=1,\ldots,m$,
and using two facts: the first one  that $\xi_j \in [a_j,b_j]$ if
and only if $\eta_j \in I_j^{\prime}$ and the second that
 \[
\sum\limits_{j=1}^m \xi_j = \sum\limits_{j=1}^m  a_j + (1-\alpha) \,
\eta_j = \alpha + 1- \alpha = 1,
 \]
we get
\begin{equation*} 
T(D)   = \Biggl\{\sum\limits_{j=1}^m \xi_j
 \, x_j \, : \ \xi_j \in [a_j,b_j], \ \sum\limits_{j=1}^m
 \xi_j = 1 \Biggr\}
  = \co(S,\cJ_S).
\end{equation*}
It follows now from Proposition~\ref{p1e} that $T(D)$ has fewer
vertices than $D$. Since $D$ has at most $n \binom{m}{n}$ vertices
we conclude that $\co(S,\cJ_S)$ has at most $n \binom{m}{n}$
vertices.

To show that the bound is best possible we provide an example in
which $\co(S,\cJ_S)$ has exactly $n \binom{m}{n}$ vertices. Take $S
= \{e_1, \ldots, e_m\}\subset \rr^m$ and
 \begin{equation*}
 \cJ_S = \bigl\{ I_1, \ldots, I_m
\bigr\}, \quad   \text{where} \quad I_j =
\left[0,\dfrac{2}{2(m-n)+1}\right], \, j=1,\ldots,m.
 \end{equation*}
Similarly as above we can check that $\co(S,\cJ_S)=D_1$ is the
intersection of $C$ and the hyperplane
\begin{equation*} 
H=\left\{ (\zeta_1,\ldots, \zeta_m) \in \nR^m\, : \, \zeta_1 +
\cdots + \zeta_m = m-n + \dfrac{1}{2}  \right\}.
\end{equation*}
One can immediately check that $H$ intersects exactly $n \binom{m}{n}$
edges of $C$ at the points whose $m-n$ coordinates are equal to $1$,
$n-1$ coordinates are equal to $0$ and exactly one coordinate is
equal to $1/2$.  It is easy to see that no three such points are
collinear.
  Similarly as in the case of $D$,  all the vertices of
$D_1$ must be the points of intersection of $H$ and edges of $C$.
 Clearly, no edge of $C$ lies on $H$.
   Therefore the intersection
of $H$ and $C$ has exactly $n \binom{m}{n}$ vertices. Since
$\co(S,\cJ_S)$ coincides with the intersection of the hyperplane $H$
and $C$ it has exactly $n \binom{m}{n}$ vertices. The proof of the
theorem is complete.
\end{proof}

\begin{remark}
In the proof of Theorem 5 in \cite{az}  the hyperplane
\begin{equation*} 
\Pi_n=\left\{ (\zeta_1,\ldots, \zeta_m) \in \nR^m \, : \, \zeta_1 +
\cdots + \zeta_m= n \right\}
\end{equation*}
is mentioned as one which makes the bound  $n \binom{m}{n}$ best
possible. In fact, this hyperplane contains $\binom{m}{n}$ vertices
of $C$. Since each vertex belongs to exactly $n$ edges it could be
argued that the hyperplane $\Pi_n$ intersects $n \binom{m}{n}$ edges
of $C$. Note that our hyperplane $H$ actually intersects $n
\binom{m}{n}$ edges of $C$ at distinct points.
\end{remark}

The subsequent theorem shows that in special cases the number of
vertices of $\co(S,\cJ_S)$ cannot be too large. We will need an
additional definition. A family $\cJ_S=\{I_j=[a_j,b_j]\;,\;\;\jom\}$
is called {\it wide}
 if $ d_k+d_j >1-\alpha $ for all $k\not=j$, where $d_i=b_i-a_i$ and
$\alpha=\sum_{i=1}^m a_i$. One can easily check that the family
$\cJ_S$ considered in the example finishing the proof of
Theorem~\ref{tmax} is wide only when $n=3$ or $n=4$ and in both
cases the maximal number of vertices of $\co(S,\cJ_S)$ guaranteed by
Theorem~\ref{tmax} is the same as the one guaranteed by the
following theorem.

\begin{theorem}\label{t9}
Let $S = \{x_1,\ldots,x_m\}$ be a subset of distinct points in
$\,\bL$. Assume that $\cJ_S=\{I_j=[a_j,b_j]\,:\,j=1,\dots,m\}$ is a
wide family of intervals,  with $a_j<b_j<1-\alpha+a_j $ and $\alpha<
1$. Then $\co(S,\cJ_S)$ is the convex hull of at most $m(m-1)$
points and this bound is best possible.
\end{theorem}

\begin{proof}
Similarly as in the proof of Theorem~\ref{tmax} we shall show that
$\co(S,\cJ_S)$ is an image upon a linear transformation of a
polytope, call it $B$,  having
 $m(m-1)$ vertices and lying in the hyperplane
\begin{equation*} 
\Pi_1 = \bigl\{ (\xi_1,\ldots, \xi_m) \in \nR^m \, : \,\,
 \xi_1+ \xi_2+\cdots +\xi_m = 1 \bigr\}.
\end{equation*}
To  construct $B$, we first consider the points
 \[
 v_j=(a_1,\dots,1-\alpha+a_j,\dots,a_m), \quad \jom,
 \]
lying on $\Pi_1$. Clearly, $\Delta=\conv\{v_1,\dots,v_m\}$ is a
fully dimensional simplex in $\Pi_1$. Define
\[
B=\Delta\cap H_1^+\cap H_2^+\cap \dots \cap H_m^+
\]
in which
 \[
 H_j^+=\{(\xi_1,\ldots, \xi_m) \in \nR^m \, : \,
  \xi_j \leq b_j\},\quad \jom,
 \]
is the halfspace bounded by the hyperplane
\[
H_j= \bigl\{(\xi_1,\ldots,\xi_m) \in \nR^m \, : \,
  \xi_j = b_j \bigr\}.
\]

There are $m-1$ edges of $\Delta$ emanating from  a   vertex $v_j$
of $\Delta$. Clearly, each one of these edges  intersects  the
hyperplane $H_j$ at a point $v_j^k$, $k\not=j$. It is easy to check
that
 \[
v_j^k=(a_1,\dots,b_j,\dots,c_k,\dots,a_m),\quad k\not=j,
 \]
where  $c_k=1-\alpha-d_j+a_k$. Thus, each intersection
$\Delta_j=\Delta\cap H_j$, $\jom$, is a simplex with $m-1$ vertices
$v_j^1,v_j^2,\dots,v_j^{j-1},v_j^{j+1},\dots,v_j^m$.

Now we shall check that  every vertex of any simplex $\Delta_j$,
$\jom$, is  a vertex of $B$. Indeed, take
$v_j^k=(a_1,\dots,b_j,\dots,c_k,\dots,a_m)$, $k\not=j$. Obviously
$v_j^k\in \Delta \cap H_t^+$ for  $t\not=k$. To show that also
$v_j^k\in \Delta\cap H_k^+$ we need to check  that $c_k$ (the $k$-th
coordinate of  $v_j^k$) is less than $b_k$. This is true since the
inequality
\[
c_k=1-\alpha-d_j+a_k<b_k
\]
is equivalent to
\[
1-\alpha<b_k-a_k+d_j=d_k+d_ j
\]
and the latter inequality is  true because the family  $\cJ_S$ is
wide. In this way we have shown that every vertex of $\Delta$ gives
rise to $m-1$ vertices of $B$. Thus $B$ is a polytope with $m(m-1)$
vertices.

Now consider  a linear transformation $T_S: \rr^m\to \bL$ defined by
\[
T_S(\xi_1,\xi_2,\dots,\xi_m) := \su\xi_j x_j.
\]
 We want to show that
\[
\co(S,\cJ_S)=T_S(B).
\]
The inclusion $T_S(B)\subseteq \co(S,\cJ_S)$ simply follows from the
definitions given above.

To show the reverse inclusion suppose to the contrary that there
exists
 \[
 z\in \co(S,\cJ_S)\setminus T_S(B).
  \]
Of course, $z=\su\mu_j x_j$ for some $\mu_1,\dots,\mu_m$ such that
$a_j\le \mu_j\le b_j$, $\jom$, and $\su\mu_j=1$. Obviously,
 \[
 (\mu_1,\dots,\mu_m)\in \Pi_1\cap\bigcap_{j=1}^m H_j^+\setminus B.
 \]
From the definition of $B$ it follows now that
$(\mu_1,\dots,\mu_m)\not\in \Delta$. As $\Delta$ is a fully
dimensional simplex in $\Pi_1$ we have
\begin{eqnarray} \label{mu1}
(\mu_1,\dots,\mu_m)\in \aff \{v_1,\dots,v_m\}\setminus \conv
\{v_1,\dots,v_m\}.
\end{eqnarray}
From (\ref{mu1}) we get
\begin{eqnarray} \label{mu2}
(\mu_1,\dots,\mu_m)=\su\lambda_jv_j
\end{eqnarray}
for some numbers $\lambda_1,\dots,\lambda_m$, satisfying
$\su\lambda_j=1$, among which at least one does not belong to the
interval $[0,1]$. In connection with the last observation we infer
that there exists   $i_0$ such that $\lambda_{i_0}<0$.  It is easy
to check that (\ref{mu2}) is equivalent to
\begin{eqnarray*}
(\mu_1,\dots,\mu_m) =
(a_1+\lambda_1(1-\alpha),\dots,a_m+\lambda_m(1-\alpha)),
\end{eqnarray*}
which gives
\begin{eqnarray*}
\mu_{i_0}=a_{i_0}+\lambda_{i_0}(1-\alpha)<a_{i_0}
\end{eqnarray*}
and  contradicts the condition $a_{i_0}\le \mu_{i_0}\le b_{i_0}$.
Thus $\co(S,\cJ_S)\subseteq T_S(B)$ and consequently
$\co(S,\cJ_S)=T_S(B)$. Therefore $\co(S,\cJ_S)$ cannot have more
than $m(m-1)$ vertices.

Next we show that the number $m(m-1)$ is attained for wide families.
Let $e_1, \ldots, e_m$ be the unit vectors in $\nR^m$.  Define
 \begin{align*}
S = & \{e_1, \ldots, e_m\}, \quad \cJ_S = \bigl\{ I_1, \ldots, I_m
\bigr\}, \\  & \text{where} \quad I_j = \bigl[0, 2/3\bigr], \,
j=1,\ldots,m.
 \end{align*}
Clearly $\cJ_S$ is a wide family and $\co(S,\cJ_S)$ has exactly
$m(m-1)$ vertices.
\end{proof}

\section{Minimal families of intervals} \label{smi}

The convex interval hull of a set $S$ essentially depends on the
family of intervals $\cJ_S$ associated with $S$. In
Example~\ref{ex1} we saw that the convex interval hull produced by
the family $\cJ_S$ of intervals $[0,t_j]$ with $t_j>1$ produces the
same convex interval hull as the family of intervals $[0,1]$.  This
observation indicates that the latter family is in some sense
minimal. In this section we define and explore the minimality of
families of intervals.

\begin{definition} \label{dmin}
Let $S$ be a finite set of points in $\bL$. A family of intervals
$\cJ_S =\bigl\{I_1, \ldots, I_m\bigr\}$ is a {\em minimal interval
family for the set} $S$ if
\[
\cJ_S^{\prime} = \bigl\{I_1^{\prime}, \ldots, I_m^{\prime}\bigr\}, \
\ I_j^{\prime} \subseteq I_j, \ j=1,\ldots,m,
\]
and
\[
\co(S,\cJ_S^{\prime}) = \co(S,\cJ_S)
\]
imply that
\[
I_j^{\prime} = I_j, \ j=1,\ldots,m.
\]
\end{definition}

\begin{definition} \label{dire}
Let $\cJ= \bigl\{ I_1,\ldots,I_m \bigr\}$ be a family of bounded
intervals such that $I_j = [a_j,b_j], \, a_j \leq b_j, \,
j=1,\ldots,m$. Set $\alpha = \su a_j$ and $\beta = \su b_j$.  The
family $\cJ$ is called {\em irreducible} if $b_k - a_k \leq
\min\{1-\alpha,\beta - 1\}$ for all $k=1,\ldots,m$.
\end{definition}

Let, as before,
\[
\cJ = \bigl\{ I_1,\ldots,I_m \bigr\}, \ \ I_j = [a_j,b_j], \ j =
1,\ldots,m, \ \ \alpha = \su a_j,  \ \ \beta = \su b_j,
\]
and assume $\alpha \leq 1 \leq \beta$. In the rest of this section
we will use the following notation. With the family $\cJ$ we
associate the following family $\wh \cJ$:
\[
\wh \cJ = \bigl\{ \wh{I}_1,\ldots, \wh{I}_m \bigr\}, \ \ \wh{I}_j =
[\wh{a}_j, \wh{b}_j], \ j=1,\ldots,m, \ \ \wh{\alpha} = \su
\wh{a}_j, \ \ \wh{\beta} = \su \wh{b}_j,
\]
where
\[
\wh{a}_j = \max\bigl\{a_j,b_j-(\beta-1)\bigr\}, \ \ \ \wh{b}_j =
\min\bigl\{b_j,a_j +(1 - \alpha)\bigr\}.
\]
Since we assume $\alpha \leq 1 \leq \beta$, we clearly have
\begin{equation} \label{eqwhs}
a_j \leq \wh{a}_j \leq \wh{b}_j \leq b_j, \ \ \ \ j=1,\ldots,m.
\end{equation}
The following implication is straightforward: if $\cJ$ is
irreducible, then $\cJ = \wh{\cJ}$.  In the next lemma we study the
relationship between $\cJ$ and $\wh{\cJ}$ further.  Among other
statements we prove the converse of the last implication.  We set
\[
\Pi_1 = \bigl\{\bigl(\zeta_1,\ldots,\zeta_m\bigr) \in \nR^m \, : \,
\zeta_1 + \cdots +\zeta_m = 1 \bigr\}.
\]

\begin{lemma} \label{lpro}
Let $\cJ = \bigl\{I_1, \ldots, I_m\bigr\}, \, I_j=[a_j,b_j],
j=1,\ldots,m$, be a family of bounded intervals. Set $\alpha =
\sum_{j=1}^m a_j$ and $\beta = \sum_{j=1}^m b_j$ and assume $\alpha
\leq 1 \leq \beta$. The following three statements hold.
\begin{enumerate}[{\rm (a)}]
\item \label{iproa}
Let $k \in \{1,\ldots,m\}$.  The projection of the set
\[
\bigl( I_1 \times \cdots \times I_m \bigr) \cap \Pi_1 \subset \nR^m
\]
to the $k$-th coordinate axes in $\nR^m$ is the interval $\wh{I}_k =
\bigl[ \wh{a}_k, \, \wh{b}_k \bigr]$.
\item \label{iprob}
$\displaystyle \bigl( \wh{I}_1 \times \cdots \times \wh{I}_m \bigr)
\cap \Pi_1 = \bigl( I_1 \times \cdots \times I_m \bigr) \cap \Pi_1$.
\item \label{iproc}
The family $\wh{\cJ}$ is irreducible.
\end{enumerate}
\end{lemma}
\begin{proof}
The statement (\ref{iproa}) claims the equality of two sets. To
prove it, let
\[
\bigl(\xi_1,\ldots,\xi_m\bigr) \in \bigl( I_1 \times \cdots \times
I_m \bigr) \cap \Pi_1.
\]
Then
\begin{equation*}
\xi_k = 1 - \sum_{\substack{j=1,j\neq k}}^m \xi_j \leq 1 -
\sum_{\substack{j=1,j\neq k}}^m a_j = 1 - \bigl(\alpha - a_k\bigr)
\end{equation*}
and
\begin{equation*}
\xi_k = 1 - \sum_{\substack{j=1,j\neq k}}^m \xi_j \geq 1 -
\sum_{\substack{j=1,j\neq k}}^m b_j = 1 - \bigl(\beta - b_k\bigr).
\end{equation*}
Since $a_k \leq \xi_k \leq b_k$, it follows that $\wh{a}_k \leq
\xi_k \leq \wh{b}_k$. This proves that the projection onto $k$-th
coordinate is contained in the interval $\wh{I}_k$.

For simplicity of notation, we will prove the converse inclusion for
$k = 1$.  Let $\xi_1\in \wh{I}_1$.  Then
\[
1- \min\bigl\{b_1,a_1 + 1 - \alpha\bigr\} \leq 1 - \xi_1 \leq 1-
\max\bigl\{a_1,b_1 - \beta + 1\bigr\}
\]
and consequently
\begin{equation} \label{eqxi1}
\alpha - a_1  \leq 1 - \xi_1 \leq \beta - b_1.
\end{equation}
Since the function
\[
\bigl(\zeta_2,\ldots,\zeta_m\bigr) \mapsto \sum_{j=2}^m \zeta_j
\]
is a continuous function on $I_2 \times \cdots \times I_m$ with the
minimum $\alpha - a_1$ and the maximum $\beta - b_1$, its range is
$\bigl[\alpha - a_1, \beta - b_1\bigr]$.  Now \eqref{eqxi1} implies
that there exists
\[
\bigl(\xi_2,\ldots,\xi_m\bigr) \in I_2 \times \cdots \times I_m
\]
such that $\sum_{j=2}^m \xi_j = 1- \xi_1$. Thus
\[
\bigl(\xi_1,\ldots,\xi_m\bigr) \in \bigl( I_1 \times \cdots \times
I_m \bigr) \cap \Pi_1,
\]
and (\ref{iproa}) is proved.

The statement (\ref{iprob}) follows from the fact that (\ref{iproa})
holds for all $k=1,\ldots,m$, and $\wh{I}_k \subseteq I_k$.

To prove (\ref{iproc}), we notice that (\ref{iprob}) implies that
for each $k=1,\ldots,m$, the projection of  $\bigl( \wh{I}_1 \times
\cdots \times \wh{I}_m \bigr) \cap \Pi_1$  to the $k$-th coordinate
axes in $\nR^m$ is the interval $\wh{I}_k = \bigl[ \wh{a}_k, \,
\wh{b}_k \bigr]$. Furthermore, an application of (\ref{iproa}) to
the family $\wh{\cJ}$ yields that the same projection is the
interval
\[
\Bigl[  \max\bigl\{\wh{a}_k,\wh{b}_k-(\wh{\beta}-1)\bigr\},
\min\bigl\{\wh{b}_k,\wh{a}_k +(1 - \wh{\alpha})\bigr\} \Bigr].
\]
Consequently,
\[
\wh{a}_k = \max\bigl\{\wh{a}_k,\wh{b}_k-(\wh{\beta}-1)\bigr\}, \ \ \
\wh{b}_k = \min\bigl\{\wh{b}_k,\wh{a}_k +(1 - \wh{\alpha})\bigr\},
\]
and hence
\[
\wh{a}_k \geq \wh{b}_k-(\wh{\beta}-1), \ \ \ \ \wh{b}_k \leq
\wh{a}_k +(1 - \wh{\alpha}).
\]
This implies that $\wh{\cJ}$ is irreducible and the lemma is proved.

\end{proof}

\begin{proposition} \label{pminJS}
Let $S$ be a finite subset of $\bL$ and let
 $\cJ_S$ be a corresponding family of bounded intervals such that
$\co(S,\cJ_S) \neq \emptyset$. Then
\[
\co(S, \wh{\cJ}_S) = \co(S,\cJ_S).
\]
\end{proposition}
\begin{proof}
The proposition follows from (\ref{iprob}) in Lemma~\ref{lpro}.
\end{proof}
\begin{theorem} \label{tmiir}
Let $S$ be a finite subset of $\bL$ and let
 $\cJ_S$ be a corresponding family of bounded intervals such that
$\co(S,\cJ_S) \neq \emptyset$. If $\cJ_S$ is a minimal family for
$S$, then $\cJ_S$ is irreducible.
\end{theorem}
\begin{proof}
We prove the contrapositive. Assume that $S$ has $m$ elements and
let $\cJ_S = \bigl\{I_1, \ldots, I_m\bigr\}$.  Assume further that
$\cJ_S$ is not irreducible. Then there exists $k \in \{1,\ldots,m\}$
such that $\wh{I}_k$ is a proper subset of $I_k$. But $\co(S,
\wh{\cJ_S}) = \co(S,\cJ_S)$ by Proposition~\ref{pminJS}. Since
$\wh{I}_j \subseteq I_j$ for all $j=1,\ldots,m$, \ $\cJ_S$ is not a
minimal family for $S$.
\end{proof}

The following example shows that the converse of Theorem~\ref{tmiir}
is not true.
\begin{example}
Consider the set $S$ of $4$ points in $\nR^2$ from
Example~\ref{ex3}. Let $\cJ_S$ be the family of four copies of the
interval $[0,1]$. Then $\co(S,\cJ_S) = \conv S$.  The family $\cJ_S$
is clearly irreducible, but it is not a minimal family for $S$,
since the family
 $
\cJ_S' = \bigl\{[0,1], [0,1], [0,1], [0,t]\bigr\},
 $
for arbitrary $t \in [0,1)$, clearly produces the same convex
interval hull.
\end{example}

In the next theorem we show that for affinely independent sets the
converse of Theorem~\ref{tmiir} holds true.  Recall that a set
$S=\{y_1,\dots,y_m\}$ of points in $\bL$ is {\em affinely
independent} if and only if the affine mapping
\begin{equation} \label{eqbij}
\bigl(\xi_1,\ldots,\xi_m\bigr) \mapsto \xi_1\,y_1+\cdots +\xi_m \,
y_m
\end{equation}
is a bijection between $\Pi_1^m$ and $\aff S$.

\begin{theorem}
Let $S$ be a finite affinely independent subset of $\bL$ and let
$\cJ_S$ be an associated family of bounded intervals. The family
$\cJ_S$ is a minimal family for $S$ if and only if it is
irreducible.
\end{theorem}

\begin{proof}
Let $S$ be an affinely independent set with $m$ elements and
$\cJ_S=\bigl\{I_1,\ldots,I_m\bigr\}$.  We prove the contrapositive
of the ``if'' part of the theorem. Assume that $\cJ_S$ is not a
minimal family for $S$.  Then there exist a family of intervals
\begin{gather} \nonumber
\cJ_S^{{\prime}}=\bigl\{I_1^{{\prime}},\ldots,I_m^{{\prime}}\bigr\}
\ \ \ \text{such that} \ \ \ I_j^{{\prime}} \subseteq I_j, \
j=1,\ldots,m, \\
 \label{eqeqcos}
 \co(S,\cJ_S^{\prime})= \co(S,\cJ_S),
\end{gather}
and there exists $k\in \{1,\ldots,m\}$ for which $I_k^{\prime}$ is a
proper subset of $I_k$. Setting $I_k^{\prime} =
\bigl[a_k^\prime,b_k^\prime\bigr],  \, I_k = \bigl[a_k,b_k\bigr]$,
the last condition is equivalent to
\begin{equation} \label{eqinIpI}
b_k^\prime - a_k^\prime < b_k - a_k.
\end{equation}
Since the mapping \eqref{eqbij} is a bijection, \eqref{eqeqcos} is
equivalent to
\begin{equation*}
\bigl(I_1^{{\prime}}\times \cdots \times I_m^{{\prime}}\bigr) \cap
\Pi_1 = \bigl(I_1\times \cdots \times I_m\bigr) \cap \Pi_1.
\end{equation*}
By Lemma~\ref{lpro}(\ref{iproa}), the last equality implies that
$\wh{\cJ^\prime_S} = \wh{\cJ_S}$. Therefore, by \eqref{eqwhs} and
\eqref{eqinIpI},
\[
\wh{b}_k - \wh{a}_k = \wh{b^\prime}_k - \wh{a^\prime}_k \leq
b_k^\prime - a_k^\prime < b_k - a_k.
\]
Hence $\wh{\cJ_S} \neq \cJ_S$, and consequently $\cJ_S$ is not
irreducible.
\end{proof}

\section{The convex interval hull and the homothety} \label{scoh}

Let $\delta$ be a nonzero real number and $v \in \bL$. The
transformation $H_v^{\delta}: \bL \to \bL$ defined by
 \[
H_v^{\delta}(x) := v+\delta \, x, \quad x \in \bL,
 \]
is called a {\it homothety}. If $\delta > 0$ the homothety is called
{\em positive} and if $\delta < 0$ the homothety is called {\em
negative}.  The number $\delta$  is called the {\it ratio of the
homothety}. The image of $K \subset \bL$ under $H_v^{\delta}$ is
denoted by $H_v^{\delta}(K)$ and it is called a {\it homothet of}
$K$. It is convenient to set $H_v^{0}(K) = \emptyset$.

Let $S = \{x_1,\ldots,x_m\}$, be a finite set of points in $\bL$. We
are interested only in homotheties that map $\aff S$ to $\aff S$.
Let $v \in \bL$ and $\delta \neq 0$. Clearly $H_v^{\delta}(\aff S)
\subseteq \aff S$ if and only if there exist $\nu_j \in \nR,\,
j=1,\ldots,m$, such that
\begin{equation} \label{eql1}
v = \su \nu_j x_j \qquad \text{and} \qquad \delta = 1 - \su \nu_j.
\end{equation}

\begin{theorem} \label{t5}
Let $S = \{x_1,\ldots,x_m\}$ be a finite set of points in $\bL$ and
let $\cJ_S = \bigl\{I_1,\ldots,I_m\bigr\}$ be a corresponding family
of nonempty intervals.  Let $v \in \bL$ and $\delta \neq 0$ be such
that $H_v^{\delta}(\aff S) \subseteq \aff S$. Assume that
\eqref{eql1} holds and set $h_j(t) = \nu_j + \delta\, t, \, t\in
\nR$. Then
\begin{equation*} 
H_v^{\delta}\bigl(\co(S, \cJ_S)\bigr) = \co\bigl(S, \cJ'_S\bigr),
\end{equation*}
where
\begin{equation*}  
\cJ'_S  = \{ I'_1, \ldots,I'_m\},\ \ I'_j = h_j\bigl(I_j\bigr), \ \
j = 1,\ldots,m.
\end{equation*}
\end{theorem}

\begin{proof}
To prove the inclusion $\co\bigl(S, \cJ'_S\bigr) \subseteq
H_v^{\delta}\bigl(\co(S, \cJ_S)\bigr)$, let $y \in \co\bigl(S,
\cJ'_S\bigr)$. Then there exist $\xi_j \in I_j,\jom$, such that
\[
y = \sum_{j=1}^m h_j\bigl(\xi_j\bigr) \, x_j \ \ \ \ \text{and} \ \
\ \ \sum_{j=1}^m h_j\bigl(\xi_j\bigr) = 1.
\]
Since, by \eqref{eql1},
\[
\sum_{j=1}^m \xi_j = \frac{1}{\delta}\Biggl(1 - \sum_{j=1}^m
\nu_j\Biggr) = 1,
\]
with $x = \sum_{j=1}^m  \xi_j\, x_j \in \co\bigl(S,\cJ_S\bigr)$ we
have
\begin{equation*}
y = \sum_{j=1}^m  \bigl(\nu_j + \delta \, \xi_j \bigr) \, x_j =
\sum_{j=1}^m  \nu_j\, x_j + \delta \sum_{j=1}^m \xi_j \, x_j  =
H_v^\delta(x).
\end{equation*}

The converse inclusion is proved similarly and the theorem is
established.
\end{proof}

\begin{corollary} \label{cchS}
Let $S = \{x_1,\ldots,x_m\}$ be a finite set of points in $\bL$ and
let $c_j \in \nR$ be such that $\gamma = \su c_j \neq 1$. Set
$h_j(t) = c_j + (1-\gamma)t$ and
\[
\cJ'_S = \bigl\{I'_1,\ldots,I'_m\bigr\} \ \ \  \text{with} \ \ \
I'_j = h_j\bigl([0,1]\bigr), \ \ \ \jom.
\]
Then
\[
\co\bigl(S,\cJ'_S\bigr) =  H^{1-\gamma}_v\bigl(\conv S \bigr), \ \ \
\text{where} \ \ \ v =  \sum_{j=1}^m c_j x_j.
\]
\end{corollary}

\begin{remark} \label{r2}
We continue to use the notation of Corollary~\ref{cchS}. Further, we
assume that $c_j \geq 0, \jom$, and $0 < \gamma < 1$. Simple algebra
yields
\[
H^{1-\gamma}_v(x) = v + (1-\gamma)x = \frac{1}{\gamma}v +
(1-\gamma)\Bigl(x- \frac{1}{\gamma} v \Bigr), \ \ \ x \in \aff S.
\]
This expression shows that $\frac{1}{\gamma}v$ is a fixed point of
$H^{1-\gamma}_v$.  Since $0 < 1- \gamma < 1$ and $\frac{1}{\gamma}v
\in \conv S$ the homotet $H^{1-\gamma}_v\bigl(\conv S\bigr)$ is a
contraction of $\conv S$ and it is completely contained in $\conv
S$.

The Gauss-Lucas theorem states that all the roots of the derivative
of a complex non-constant polynomial $p$ lie in the convex hull of
the roots of $p$, called the Lucas polygon of $p$. The reasoning
presented in Remark~\ref{r2} was used in \cite{b} to improve the
Gauss-Lucas theorem by proving that all the nontrivial roots of the
derivative of $p$ lie in a convex polygon that is a strict
contraction of the Lucas polygon of $p$ and that is completely
contained in it.
\end{remark}

We conclude this article with a result  motivated by
Examples~\ref{ex22}, \ref{ex23}, \ref{ex28} and \ref{ex29}.
It is
clear that the convex interval hull  $\co(S,\cJ_S)$ in Figure~\ref{f2} is the
closure of a set difference of
$\conv S$ and the union of two smaller homotets of $\conv S$.
Similarly, $\co(S,\cJ_S)$ in Figure~\ref{f1} is the closure of a set
difference of a large homotet of $\conv S$ and the union of two
smaller homotets of $\conv S$.  The reader will easily observe
analogous properties of the convex interval hulls in  Figures~\ref{f3}, \ref{f13} and \ref{f14}.
In Theorem~\ref{cllin} below we give a general result which explains
these observations.

\begin{lemma} \label{lbc}
Let $\cJ = \bigl\{I_1, \ldots, I_m\bigr\}, \, I_j=[a_j,b_j],
j=1,\ldots,m$, be an irreducible family of intervals. Set $\alpha =
\sum_{j=1}^m a_j$ and $\beta = \sum_{j=1}^m b_j$ and assume $\alpha
\leq 1 \leq \beta$. For $k, j \in \{1,\ldots,m\}$ define
\begin{align*}
I^0_j & = \bigl[a_j,a_j + (1-\alpha)\bigr],  \\[8pt]
I^k_j & = \begin{cases} \bigl[a_j,a_j + (1-\alpha) -(b_k-a_k) \bigr]
\ \ & \text{for} \ \ \  j \neq k, \\[6pt]
 \bigl(b_j,a_j + (1-\alpha)\bigr] & \text{for} \ \ \  j = k.
\end{cases}
\end{align*}
Set
\begin{equation*}
B = \bigl(I_1\times \cdots \times I_m \bigr) \bigcap \Pi_1, \ \ \
B_u =
 \Pi_1 \bigcap \, \bigcup_{k=1}^m \bigl(I_1^k\times \cdots
\times I_m^k \bigr).
\end{equation*}
Then $B \cap B_u = \emptyset$ and
\begin{equation} \label{eqseq}
B \cup B_u  = \bigl(I_1^0\times \cdots \times I_m^0 \bigr) \bigcap
\Pi_1.
\end{equation}
\end{lemma}
\begin{proof}
The equality $B \cap B_u = \emptyset$ is obvious. Since $\cJ$ is
irreducible we have $b_j \leq a_j + 1-\alpha$ for all
$j=1,\ldots,m$. Consequently $B\cup B_u$ is a subset of
$\bigl(I_1^0\times \cdots \times I_m^0 \bigr) \cap \Pi_1$.

To prove the converse inclusion in \eqref{eqseq}, let
\begin{equation} \label{eqin0}
(\xi_1,\ldots,\xi_m) \in \bigl(I_1^0\times \cdots \times I_m^0
\bigr) \cap \Pi_1
\end{equation}
and assume that
\begin{equation} \label{eqinI}
(\xi_1,\ldots,\xi_m) \notin \bigl(I_1\times \cdots \times I_m \bigr)
\cap \Pi_1.
\end{equation}
Then there exists $k \in \{1,\ldots,m\}$ such that
\begin{equation} \label{eqke}
\xi_k \in \bigl(b_k,a_k+(1-\alpha)\bigr].
\end{equation}
Next we prove the implication
\begin{equation} \label{eqimp}
\sum_{j=1}^m \xi_j = 1 \ \Rightarrow \ \xi_j \leq a_j+(1-\alpha)
-(b_k-a_k) \ \Bigl(\!\forall\, j \in \{1,\ldots,m\} \setminus
\{k\}\! \Bigr).
\end{equation}
Since the contrapositive is easier to prove, assume
\begin{equation} \label{eqle}
\exists\, l \in \{1,\ldots,m\} \setminus \{k\} \ \ \ \text{such
that} \ \  \xi_l > a_l+(1-\alpha) -(b_k-a_k).
\end{equation}
Then, using \eqref{eqke} and \eqref{eqle}, we find
\begin{equation*}
\sum_{j=1}^m \xi_j > b_k + \bigl(a_l + (1-\alpha) -(b_k-a_k) \bigr)
+ \bigl(\alpha - a_k - a_l \bigr) = 1,
\end{equation*}
and \eqref{eqimp} is proved. Hence, we have shown that
\eqref{eqin0}, \eqref{eqinI} and \eqref{eqke} imply that
\[
\xi_j \in \bigl[a_j, a_j+ \leq (1-\alpha) -(b_k-a_k)\bigr] \ \ \
\Bigl(\!\forall\, j \in \{1,\ldots,m\} \setminus \{k\}\! \Bigr).
\]
This together with \eqref{eqke} implies that $(\xi_1,\ldots,\xi_m)
\in B_u$ and the lemma is proved.
\end{proof}

\begin{theorem} \label{cllin}
Let  $S=\{x_1,...,x_m\}$ be an affinely independent set in $\bL$.
Let $\cJ_S = \bigl\{I_1, \ldots, I_m\bigr\}, \, I_j=[a_j,b_j],
j=1,\ldots,m$, be an irreducible family of intervals and assume
$\alpha = \sum_{j=1}^m a_j < 1$. Then $\co\bigl(S,\cJ_S\bigr)$ is
the closure of the set difference of the sets
\begin{equation*}
H_v^{\delta}\bigl(\conv S\bigr) \ \ \ \text{and} \ \ \
\bigcup_{j=1}^m  H_{v+d_jx_j}^{\delta-d_j}\!\bigl(\conv S\bigr)
\end{equation*}
where $v=\su a_jx_j,\, \delta = 1 - \alpha$, and $d_j=b_j-a_j,\,
j=1,\ldots,m$.
\end{theorem}
\begin{proof}
The claim of the theorem is equivalent to the equality
\begin{equation} \label{obc}
\co\bigl(S,\cJ_S\bigr) \bigcup \ \bigcup_{j=1}^m
 H_{v+d_jx_j}^{\delta-d_j}\!\bigl(\conv S\bigr) =
H_v^{\delta}\bigl(\conv S\bigr)
\end{equation}
together with the condition that the set $\co\bigl(S,\cJ_S\bigr)$
has no common interior points with the polytopes
$H_{v+d_jx_j}^{\delta-d_j}\!\bigl(\conv S\bigr),\, \jom$. To prove
\eqref{obc} and the stated condition we use Lemma~\ref{lbc} and the
notation introduced there. For $k = 1,\ldots, m$, set
\begin{align*}
\cJ^0_S  = \bigl\{I^0_1, \ldots, I^0_m \bigr\}, \ \ \ \ \cJ^k_S  =
\bigl\{I^k_1, \ldots, I^k_m \bigr\}, \ \ \ \ \overline{\cJ}^k_S  =
\bigl\{\overline{I}^k_1, \ldots, \overline{I}^k_m \bigr\}.
\end{align*}
Since $S$ is affinely independent the affine mapping
\begin{equation}  \label{eqbij2}
\bigl(\xi_1,\ldots,\xi_m\bigr) \mapsto \xi_1\,y_1+\cdots +\xi_m \,
y_m
\end{equation}
is a bijection between $\Pi_1$ and $\aff S$. Together with
Lemma~\ref{lbc} this implies
\begin{equation} \label{eqJs}
\co\bigl(S,\cJ_S\bigr)\; \bigcup \ \bigcup_{k=1}^m
\co\bigl(S,\cJ^k_S\bigr) = \co\bigl(S,\cJ^0_S\bigr)
\end{equation}
and, for $k=1,\ldots,m$,
\[
\co\bigl(S,\cJ_S\bigr) \bigcap \co\bigl(S,\cJ^k_S\bigr) = \emptyset.
\]
Since \eqref{eqbij2} defines a continuous mapping it follows that
$\co\bigl(S,\overline{\cJ}^k_S\bigr)$ is a closure of
$\co\bigl(S,\cJ^k_S\bigr)$. Therefore, for $k=1,\ldots,m$, the
polytopes $\co\bigl(S,\overline{\cJ}^k_S\bigr)$ and
$\co\bigl(S,\cJ_S\bigr)$ have no common interior points.

By Corollary~\ref{cchS}
\begin{equation*}
\co\bigl(S,\overline{\cJ}^k_S\bigr)= H_{v+d_k
x_k}^{\delta-d_k}\!\bigl(\conv S\bigr) \ \ \ \text{and} \ \ \
\co\bigl(S,\cJ^0_S\bigr) = H_v^{\delta}\bigl(\conv S\bigr).
\end{equation*}
Substituting the last equations in \eqref{eqJs} we get \eqref{obc}.
The theorem is proved.
\end{proof}

\end{document}